\documentclass[11pt,a4paper,reqno]{amsart}
\usepackage{amsfonts,amsthm, amscd, epsfig, amsmath, amssymb,enumerate}
\usepackage{amstext}
\usepackage{xspace}
\usepackage{graphicx}
\usepackage{url}
\usepackage{enumerate}
\usepackage[all,2cell,ps]{xy}
\usepackage{xcolor}
\usepackage{subcaption}
\graphicspath{{./figure/}}
\usepackage{cite}
\usepackage{afterpage} 
\usepackage{faktor}
\usepackage[colorinlistoftodos]{todonotes}

\usepackage{dsfont}

\newtheorem{itlemma}{Lemma}[section]
\newtheorem{itproposition}[itlemma]{Proposition}
\newtheorem{itfact}[itlemma]{Fact}
\newtheorem{theorem}[itlemma]{Theorem}
\newtheorem{itcorollary}[itlemma]{Corollary}
\newtheorem{itremark}[itlemma]{Remark}
\newtheorem{itremarks}[itlemma]{Remarks}
\newtheorem{itdefinition}[itlemma]{Definition}
\newtheorem{itexample}[itlemma]{Example}

\newenvironment{fact}{\begin{itfact}\rm}{\end{itfact}}
\newenvironment{claim}{\begin{itclaim}\rm}{\end{itclaim}}
\newenvironment{lemma}{\begin{itlemma}}{\end{itlemma}}
\newenvironment{remark}{\begin{itremark}\rm}{\end{itremark}}
\newenvironment{remarks}{\begin{itremarks} \rm}{\end{itremarks}}
\newenvironment{corollary}{\begin{itcorollary}}{\end{itcorollary}}
\newenvironment{proposition}{\begin{itproposition}}{\end{itproposition}}
\newenvironment{definition}{\begin{itdefinition}\rm}{\end{itdefinition}}
\newenvironment{example}{\begin{itexample}\rm}{\end{itexample}}
\newcommand{\be}[1]{\begin{equation}\label{#1}}
\newcommand{\ee}{\end{equation}}
\newcommand{\bl}[1]{\begin{lemma}\label{#1}}
\newcommand{\br}[1]{\begin{remark}\label{#1}}
\newcommand{\brs}[1]{\begin{remarks}\label{#1}}
\newcommand{\bt}[1]{\begin{theorem}\label{#1}}
\newcommand{\bd}[1]{\begin{definition}\label{#1}}
\newcommand{\bp}[1]{\begin{proposition}\label{#1}}
\newcommand{\bfact}[1]{\begin{fact}\label{#1}}
\newcommand{\bc}[1]{\begin{corollary}\label{#1}}
\newcommand{\bex}[1]{\begin{example}\label{#1}}
\newcommand{\ec}{\end{corollary}}
\newcommand{\efact}{\end{fact}}
\newcommand{\eex}{\end{example}}
\newcommand{\el}{\end{lemma}}
\newcommand{\er}{\end{remark}}
\newcommand{\ers}{\end{remarks}}
\newcommand{\et}{\end{theorem}}
\newcommand{\ed}{\end{definition}}
\newcommand{\ep}{\end{proposition}}
\newcommand{\epr}{\end{proof}}
\newcommand{\bpr}{\begin{proof}}
\newcommand{\bcl}[1]{\begin{claim}\label{#1}}
\newcommand{\ecl}{\end{claim}}

\newcommand{\ecs}{\end{corollary}}
\newcommand{\eers}{\end{exercise}}
\newcommand{\eexs}{\end{example}}
\newcommand{\eems}{\end{example}}
\newcommand{\els}{\end{lemma}}
\newcommand{\eles}{\end{lemmaex}}
\newcommand{\ets}{\end{theorem}}
\newcommand{\eds}{\end{definition}}
\newcommand{\eps}{\end{proposition}}

\newcommand{\bi}{\begin{itemize}}
\newcommand{\ei}{\end{itemize}}
\newcommand{\ben}{\begin{enumerate}}
\newcommand{\een}{\end{enumerate}}

\def\vbar{\mathchoice{\vrule height6.3ptdepth-.5ptwidth.8pt\kern-.8pt}
   {\vrule height6.3ptdepth-.5ptwidth.8pt\kern-.8pt}
   {\vrule height4.1ptdepth-.35ptwidth.6pt\kern-.6pt}
   {\vrule height3.1ptdepth-.25ptwidth.5pt\kern-.5pt}}
\def\fudge{\mathchoice{}{}{\mkern.5mu}{\mkern.8mu}}
\def\bbc#1#2{{\rm \mkern#2mu\vbar\mkern-#2mu#1}}
\def\bbb#1{{\rm I\mkern-3.5mu #1}}
\def\bba#1#2{{\rm #1\mkern-#2mu\fudge #1}}
\def\bb#1{{\count4=`#1 \advance\count4by-64 \ifcase\count4\or\bba A{11.5}\or
   \bbb B\or\bbc C{5}\or\bbb D\or\bbb E\or\bbb F \or\bbc G{5}\or\bbb H\or
   \bbb I\or\bbc J{3}\or\bbb K\or\bbb L \or\bbb M\or\bbb N\or\bbc O{5} \or
   \bbb P\or\bbc Q{5}\or\bbb R\or\bbc S{4.2}\or\bba T{10.5}\or\bbc U{5}\or
   \bba V{12}\or\bba W{16.5}\or\bba X{11}\or\bba Y{11.7}\or\bba Z{7.5}\fi}}

\def\1{{\bf 1}}

\def\ec{\`e }

\usepackage{mathtools}
\DeclarePairedDelimiter\ceil{\lceil}{\rceil}

\begin{document}

\title[Multiple drawing and random addition]{Urn models with random multiple drawing \\ and random addition}

\begin{abstract}
  We consider a two-color urn model with multiple drawing and random
  time-dependent addition matrix. The model is very general with
  respect to previous literature: the number of sampled balls at each
  time-step is {\em random}, the addition matrix is {\em not balanced}
  and it has general {\em random entries}. For the proportion of balls
  of a given color, we prove almost sure convergence results. In
  particular, in the case of equal reinforcement means, we prove
  fluctuation theorems (through CLTs in the sense of stable
  convergence and of almost sure conditional convergence, which are
  stronger than convergence in distribution) and we give asymptotic
  confidence intervals for the limit proportion, whose distribution is
  generally unknown.
\end{abstract}

\maketitle

    \centerline{Irene Crimaldi\footnote{IMT School for Advanced
        Studies Lucca, Piazza San Ponziano 6, 55100 Lucca, Italy,
        \url{irene.crimaldi@imtlucca.it}},
      Pierre-Yves~Louis\footnote{PAM UMR 02.102, Universit\'e
        Bourgogne Franche-Comt{\'e}, AgroSup Dijon, 1 esplanade Erasme,
        F-21000, Dijon, France,
        \url{pierre-yves.louis@agrosupdijon.fr}}\footnote{Institut de
        Math\'ematiques de Bourgogne, UMR 5584 CNRS, Universit\'e
        Bourgogne Franche-Comt{\'e}, F-21000, Dijon, France,
        \url{pierre-yves.louis@math.cnrs.fr}}, Ida
      G.~Minelli\footnote{Dipartimento di Ingegneria e Scienze
        dell'Informazione e Matematica, Universit\`a degli Studi
        dell'Aquila, Via Vetoio (Coppito 1), 67100 L'Aquila, Italy,
        \url{idagermana.minelli@univaq.it}}}

\medskip


\noindent Secondary:  62P25;  91D30   ;  92C60

\date{\today}

\tableofcontents

\section{Introduction}
\label{intro}

\emph{Reinforcement} (see~\cite{pem} for a
review) means the tendency of a stochastic evolution to increase (or
sometimes decrease, so called, negative reinforcement) the occurrence
of an event in relationship with the number of time this event took
place in the past.  The P{\'o}lya urn stochastic process is the
fundamental and paradigmatic example. It led to several
generalizations.
\smallskip

The original evolution rule of the P{\'o}lya urn is based on picking
one ball in an urn filled with colored balls and replacing that ball
in the urn together with one or more balls, according to some
''updating matrix''.  More generalized samples have been considered,
leading to \emph{multi-drawing} based updating rules. In these models,
many balls are selected at each time and returned before adding some
new ones according to a reinforcement rule. Bi-color and multi-color
models have been considered, as well as models where the extraction of
the balls is with or without replacement. The number of sampled balls
is always a fixed constant and the ``replacement matrix'' is in
general assumed to be balanced, that is, the number of added balls to
the urn is constant along time
(\textit{e.g.} \cite{chen-kuba,chen-wei-2005,Higueras2006,
  Johnson2004,kuba2016classification,KubaMahmoud-balanced-affine-2017,mailler,mahmoud_2013_multisets}). In
particular, in~\cite{Idriss2018,mailler,mahmoud_2013_multisets} the
number of added balls is a deterministic function of the composition
of the extracted sample. Results deal with the asymptotic behavior,
evolution of moments, almost sure convergence and Central Limit
Theorems (CLTs) for the fraction of balls of a given color in the
urn. In the model considered in ~\cite{launay2012urns}, $m$~balls are
sampled at a time, with replacement, and the distribution of the
increment of one color follows, given the past, a binomial
distribution with parameters~$m$ and~$p$, where $p$ depends on weights
associated to the drawn colors. Results mainly deal with regimes where
``fixation'' happens, which is more interesting for reinforced random
walks applications.  Moreover, different urn models with multi-drawing
were considered in relationship with some specific applications. See
for instance
\cite{KubaMahmoudPanholzer,Kuba-Sulzbach,laslier2017,Crimaldi2008}.
\smallskip

Other urn models merge multi-drawing and random replacement matrix.
The paper~\cite{AguechSelmi-unbalanced} is a generalization
of~\cite{AguechLasmarSelmi2019} and it deals with a constant sample
size and a random replacement matrix.  This matrix can be of P{\'o}lya
(diagonal) or Friedman (anti-diagonal, reinforcement of the non chosen
color) type and its entries have time-homogenous distribution. In
particular, we point out that CLTs are not proven for the P\'olya type
case. As we will see later on, we here fill in this gap.
\smallskip

The papers~\cite{perron, crimaldi2016} study a multi-drawing model
(called HRRU, hypergeometric randomly reinforced urn model) with a
{\em random} number $N_n$ of sampled balls and a random replacement
matrix of rank~1 (bicolor case).  The number of added balls of a given
color is proportional to the number of balls of the same color in the
sample, but the {\em random reinforcement} factor is the same for both
colors. Note that this model generalizes the one recently given
in~\cite{Chen2020}. The almost sure convergence of the color
proportions toward a non degenerate random variable is
proven. Necessary and sufficient conditions for no-atoms in the
limiting distribution are given.
\smallskip

In this paper, we consider a two-color urn model, with {\em multiple
  drawing and random time-dependent addition matrix}. The model is
very general with respect to previous literature: the number of
sampled balls at each step is random, the addition matrix, defining
the number of additional balls, has general random entries. More
precisely, for both colors, the random number of added balls is
proportional to the number of balls of the same color in the sample,
with {\em possibly different random} coefficients $A_n$, $B_n$ (which
may be {\em correlated} and their distribution may {\em depend on
  time} $n$).  The model studied in~\cite{perron,crimaldi2016}
corresponds to the particular case $A_n=B_n$.  The reinforcement rule
we consider is {\em not balanced} (thus the long-run behavior of the
total number~$S_n$ of balls in the urn at time $n$ needs to be
studied). We prove almost sure convergence results for the proportion
as well as fluctuation results, through central limit theorems in the
sense of stable convergence and of almost sure conditional
convergence, by suitably extending some approaches employed in the urn
model literature without multi-drawing (see \cite{Berti2010,
  BCPR-urne1, mayflournoy}).  Specifically, we consider two cases. If
the factors $A_n$ and $B_n$ have the same mean (equal reinforcement
means case), the limit proportion~$Z$ is random without atoms. In the
case of unequal limit reinforcement means, the proportion converges
almost surely to~$1$ (or~$0$).  When the limit proportion~$Z$ is
random, the proven central limit theorems are employed in order to
obtain asymptotic confidence intervals.
\smallskip

Some applications of the urn models with multi-drawing are described
in \cite{Kuba-Sulzbach}. Moreover, like explained
in~\cite{perron,crimaldi2016}, the present model may be applied in the
context of \emph{technology adoption} to model, for example, the
evolution of the choice between different operative systems by
companies. Below we illustrate other possible interpretations in the
contexts of \emph{opinion dynamics} and propagation of contagious
diseases (\emph{epidemic models}).
\smallskip 

Applications to opinion dynamics could be developed as follows.
Assume to be before an election between two candidates. People decide
who they are going to vote for.  People who have already decided are
represented as the colored balls already in the urn, the color meaning
the choice for one candidate. One assume this is a not evolving
choice. At each iteration, a group (with random size $N_n$) of people
is sampled (without replacement) and each one is given the opportunity
to convince a group of other people.  The new-comers will adopt the
same choice as the person who convinced them. The heterogeneity of
this reproduction mechanism is modeled through the time-dependent
randomness of the factors $A_n$ and $B_n$. The assumption of equal
reinforcement means would mean that in the long-run no advantage is
given to any party.  We can also consider the evolution of the
diffusion of a binary opinion through social networks, like \emph{Twitter}.
Each agent inside a connected community has an un-changing opinion
(for instance, a vote or a purchased product).  This community will
grow dynamically through immigration of followers.  At each step, a
subset (with random size $N_n$) of agents is chosen. Each agent of
this committee is allowed to call into the community of followers sharing
their opinion. Once again, the heterogeneity of this growth mechanism
is modeled by allowing the multiplying factors~$A_n$ and $B_n$ for
each opinion to be random. Correlation between these growth
coefficients are possible. If one of these coefficients is eventually
larger in mean, then the associated opinion will dominate eventually
(but may take some time). If both coefficients are equal in mean then
some random equilibrium takes place.
\smallskip 
 
In the original paper~\cite{egg-pol}, where the P\'olya urn model was
first defined, smallpox epidemy was the context it was applied to (see
for instance~\cite{mah,Kotz1997} and references therein).  Therefore,
a second application of our model one could have in mind is the
diffusion of genetic variants of viruses (see for
instance~\cite{EpidemicNetworks} for a review on epidemic models on
networks). We do not pretend to do any modeling study here but want to
illustrate the potentialities of our model as a ``toy model''. Assume
one want to model the propagation of a virus, existing in two
forms. Assume to consider a time scale such that there are infinitely
many persons to be possibly contaminated and that once a person is
contaminated, he/she remains contagious ``for ever'' (no recovering,
no dying).  Balls in the urn represent the contaminated persons by one
of the two variants of the virus (corresponding to the two possible
colors of the balls). We do have in mind the initial exponential
regime of the propagation of two competing variants of one virus.
Each discrete time-step of the urn's evolution means a contagion
step. People that are contaminating are assimilated to the sample made
without replacement in the urn. This is a random number $N_n$ and this
randomness may depend on time and on the total number of contaminated
persons. One chosen contaminating person diffuse the same
variant. Each variant has its own amplifying factor~$A_n$
(resp.~$B_n$): one assume that each selected person, contaminated by a
given variant, is contaminating the same number of people. This
somewhat unrealistic hypothesis is compensated by the fact that the
number of individuals infected by one person is random, with a
time-dependent and variant-dependent randomness. Moreover, $A_n$ and
$B_n$ could be correlated. This model gives insights: if the limit
means (time-asymptotic reproduction means of each variant in this
context) are unequal, one kind of virus will eventually dominate. If
they are equal, there is a limiting genuinely random proportion, for
which we provide confidence intervals.
%
%
\smallskip 

Finally, another application context could be population dynamics in
case of competitive or cooperative growth.  As before, the flexibility
of the model lies in the choice of $N_n$, $A_n$ and $B_n$. The joint
distribution of $[A_n,B_n]$ is important to model competition or
cooperation. One may think to bacterial populations and the evolution
of their respective proportions in the microbial gut.
\smallskip

The paper is organized as follows.  In Section~\ref{model} we formally
define the model.  In Section~\ref{results} we state and prove the
main results.  In Subsection~\ref{subsec:as-cv} we prove the almost
sure convergence towards a limit proportion~$Z$.  Different behaviors
occur according to equality/unequality of the limit reinforcement
means.  In particular, in the case of equal reinforcement means, we
provide precise asymptotic rates: indeed, in
Subsection~\ref{subsec:CLT} we establish central limit theorems for
the proportion~$Z_n$ of the balls of a given color in the urn and for
the empirical mean~$M_n$ of the proportion of the balls of a given
color in the samples.  Moreover, in the case of equal reinforcement
means, in Subsection~\ref{subsec:no-atom}, we prove that the
distribution of the limit proportion~$Z$ has no atoms and, in
Subsection~\ref{conf-int}, we provide asymptotic confidence intervals
for~$Z$, centered in $Z_n$ and $M_n$. We then present in Section
\ref{sec:simulations} more specific examples, illustrated with some
numerical simulations.  The paper is enriched with an appendix in
three parts which collects some more technical lemmas and general
results, in particular about stable convergence and its variants.
 

\section{The model}
\label{model}

An urn contains $a\in{\mathbb N}\setminus\{0\}$ balls of color A and
$b\in{\mathbb N}\setminus\{0\}$ balls of color B. At each discrete 
time $n\geq 1$, we simultaneously (\textit{i.e.} without replacement) draw a
random number $N_n$ of balls. Let $X_n$ be the number of extracted
balls of color A.  Then we return the extracted balls in the urn
together with other $A_nX_n$ balls of color A and $B_n(N_n-X_n)$ balls
of color B. More precisely, we take a probability space
$(\Omega,\mathcal{A},P)$ and, on it, some random variables
$N_n,\,X_n,\,A_n,\,B_n$ such that, for each $n\geq 1$, we have:

\begin{itemize}
\item[(A1)] The conditional distribution of the random variable $N_n$ given 
$$[N_1,X_1, A_1,B_1\dots,N_{n-1}, X_{n-1}, A_{n-1},B_{n-1}]$$ is
  concentrated on $\{1,\dots,S_{n-1}\}$ where $S_{n-1}$ is the total
  number of balls in the urn at time~$n-1$, that is 
\begin{equation}
S_{n-1}=a+b+\sum_{j=1}^{n-1} A_j X_j+\sum_{j=1}^{n-1}B_j(N_j-X_j).
\end{equation}

\item[(A2)] The conditional distribution of the random variable $X_n$ given 
$$[N_1,X_1, A_1,,B_1\dots,N_{n-1}, X_{n-1}, A_{n-1},B_{n-1}, N_n]$$ is
  hypergeometric with parameters $N_n,\, S_{n-1}$ and $H_{n-1}$, 
  where $H_{n-1}$ is the total number of balls of color~A at time~$n-1$, that is 
\begin{equation}
H_{n-1}=a+\sum_{j=1}^{n-1} A_j X_j.
\end{equation}

\item[(A3)] The random vector $[A_{n},B_n]$ takes values in ${\mathbb
  N}\setminus\{0\}\times{\mathbb N}\setminus\{0\}$ and it is
  independent of
  $$[N_1,X_1,A_1,B_1,\ldots,N_{n-1}, X_{n-1},A_{n-1},B_{n-1}, N_n, X_n]\,.$$
  \end{itemize}

\indent According to the above notation, the random variable $X_n$
corresponds to the number of balls having the color A in a random
sample without replacement of size $N_n$ from an urn with $H_{n-1}$
balls of color A and $K_{n-1}=(S_{n-1}-H_{n-1})$ balls of color B.
The reinforcement rule is of the ``multiplicative'' type: indeed, each
time $n$, we add to the urn $A_nX_n$ balls of color $A$ and
$B_n(N_n-X_n)$ balls of color $B$. Therefore, the total number of
added balls to the urn, that is $A_nX_n+B_n(N_n-X_n)$, is random and
depends on $n$.\\

\indent Note that we do not specify the conditional distribution of
the random variable $N_n$ (the sample size) given the past
$$[N_1,X_1, A_1,B_1\dots,N_{n-1}, X_{n-1}, A_{n-1},B_{n-1}]$$ nor the
distribution of $[A_n,B_n]$ (the random reinforcement factors $A_n$
and $B_n$ may have different distributions, they may be correlated and
their joint and marginal distributions may vary with $n$).\\

\indent It is worthwhile to remark that this model include the
Hypergeometric Randomly Reinforced Urn (HRRU) studied
in~\cite{perron,crimaldi2016} (take $A_n=B_n$ for all $n$), which in
turn include the model recently given in \cite{Chen2020}. In
particular, two special cases are the classical P\'olya urn (the case
with $N_n=1$ and $A_n=B_n=k\in \mathbb{N}\setminus\{0\}$ for each $n$)
and the $2$-colors randomly reinforced urn with the reinforcements for
the two colors equal or different in mean (the case with $N_n=1$ for
each $n$ and $[A_n, B_n]$ arbitrarily random in ${\mathbb
  N}\setminus\{0\}\times{\mathbb N}\setminus\{0\}$).  Moreover, as
told in Section \ref{intro}, previous literature (we refer to the
quoted papers in Sec.~\ref{intro}) deals with the case when the sample
size $N_n$ is a fixed constant, not depending on $n$, and/or the
balanced case (constant number of added balls to the urn each time).
\\

\indent We set $Z_n$ equal to the proportion of balls of color A in the
urn (immediately after the updating of the urn at time $n$ and
immediately before the $(n+1)$-th extraction), that is $Z_0=a/(a+b)$
and
\begin{gather*}
  Z_n=\frac{H_n}{S_n}\quad\hbox{for } n\geq 1.
\end{gather*}
Moreover we set 
\begin{equation*}
\mathcal{F}_0=\{\emptyset, \Omega \},\quad\mathcal{F}_n=
\sigma\bigl(N_1,X_1,A_1,B_1,\ldots,N_n, X_n,A_n,B_n\bigr)
\quad\hbox{for } n\geq 1\,,
\end{equation*}
and 
\begin{equation*}
\mathcal{G}_n=\mathcal{F}_n\vee\sigma(N_{n+1}),\quad
\mathcal{H}_n=\mathcal{G}_n\vee\sigma(A_{n+1},B_{n+1})\quad\hbox{for } n\geq 0.
\end{equation*}

By the above assumptions and notation, we have
\begin{equation}\label{sp-cond-rinforzi}
E[A_{n+1}\,|\,\mathcal{G}_n]=E[A_{n+1}]\,,\qquad
E[B_{n+1}\,|\,\mathcal{G}_n]=E[B_{n+1}]
\end{equation}
and
\begin{equation}\label{sp-cond}
  \begin{split}
  E[X_{n+1}\,|\,\mathcal{H}_n]&=E[X_{n+1}\,|\,\mathcal{G}_n]=N_{n+1}Z_n\,,\\
  E[N_{n+1}-X_{n+1}\,|\,\mathcal{H}_n]&=E[N_{n+1}-X_{n+1}\,|\,\mathcal{G}_n]=
  N_{n+1}(1-Z_n)\,.
  \end{split}
\end{equation}
Finally, we set $\mathcal{X}_{n}=\{0\vee N_{n}-(S_{n-1}-H_{n-1}),\dots,N_{n}\wedge
H_{n-1}\}$ and, for each $k\in \mathcal{X}_{n}$,   
\begin{equation}\label{hypergeometric-distr}
p_{n,k}=p_k(N_{n}, S_{n-1}, H_{n-1})=
\frac{\binom{H_{n-1}}{k}\binom{S_{n-1}-H_{n-1}}{N_{n}-k}}{\binom{S_{n-1}}{N_{n}}}\,.
\end{equation}

\section{Asymptotic results}
\label{results}
In this section we prove some convergence results for the model
described in Section~\ref{model} by suitably extending some approaches
employed in the urn model literature without multi-drawing (see
 \cite{Berti2010, BCPR-urne1, mayflournoy}).  \\

\indent Set $E[A_n]=m_{A,n}$ and $E[B_n]=m_{B,n}$ for all $n$.  We
will assume that the two sequences $(m_{A,n})_n$ and $(m_{B,n})_n$
respectively converge to $m_A\in (0,+\infty)$ and $m_B\in
(0,+\infty)$.  Moreover, we will consider the following cases:
\begin{itemize}
\item[1)] $m_A>m_B$.
\item[2)] $m_{A,n}=m_{B,n}=m_n$ and so $m_A=m_B=m\in (0,+\infty)$.
\end{itemize}
For simplicity, throughout the paper, we will assume
$$
A_n\vee B_n\vee N_n\leq C\qquad\mbox{ for some (integer) constant } C.
$$
We will signal when this assumption can be easily removed.
Sometimes it may be replaced by an assumption of uniformly integrability,
but we will not focus on this fact.
\\

\indent We start with proving a result valid for both cases.

\begin{lemma}\label{colori-a-infinito}
We have 
  $$
  H_n\stackrel{a.s.}\longrightarrow +\infty\qquad\mbox{and}\qquad
 K_n=(S_n-H_n)\stackrel{a.s.}\longrightarrow +\infty\,.
  $$
\end{lemma}
As a consequence, we obviously have $S_n\stackrel{a.s.}\longrightarrow
+\infty$.
\begin{proof} 
First suppose $a\wedge b\geq C$ so that $N_{i}\leq H_{i-1}$ for each
$n$.  Let $T=\inf\{n:X_n\neq N_n\}=\inf\{n: (N_n-X_n)>0\}$.  For each
$k\geq 1$, we have
  \begin{equation*}
    \begin{split}
t_k=P\{T>k\}&=P(X_i=N_i\,,i=1,\dots, k)=
E\left[\prod_{i=1}^k \frac{H_{i-1}}{S_{i-1}}\times\dots\times
\frac{H_{i-1}-(N_i-1)}{S_{i-1}-(N_i-1)}\right]\\
&=E\left[\prod_{i=1}^k \prod_{j=0}^{N_i-1}
\frac{a-j+\sum_{h=1}^{i-1}A_hN_h}{a+b-j+\sum_{h=1}^{i-1}A_hN_h}\right]\,.
    \end{split}
  \end{equation*}
  We recall that, given $c_1,c_2,c_3>0$, we have 
  $$ x\leq c_1 \Leftrightarrow \frac{c_2+x}{c_2+c_3+x}\leq
  \frac{c_1+c_2}{c_1+c_2+c_3}\,.
  $$ Therefore, applying the above inequality with
  $x=\sum_{h=1}^{i-1}A_hN_h\leq (i-1)C^2=c_1,\, c_2=a-j, \, c_3=b$,
  we get 
\begin{equation*}
\begin{split}
  &t_k\leq
  E\left[
    \prod_{i=1}^k \prod_{j=0}^{N_i-1} \frac{a-j+(i-1)C^2}{a+b-j+(i-1)C^2}
    \right]
  \leq
  E\left[ \prod_{i=1}^k \left( \frac{a+(i-1)C^2}{a+b-N_i+1+(i-1)C^2} \right)^{N_i}
\right]
    \leq
  \\
&\prod_{i=1}^k \frac{a+(i-1)C^2}{a+b-C+1+(i-1)C^2}=
\exp\left(\sum_{i=1}^k \ln(1-(b-C)/(a+b-C+1+(i-1)C^2))\right)
\longrightarrow 0\quad\mbox{as } k\to +\infty\,.
\end{split}
\end{equation*}
This fact means that $P(T=+\infty)=\lim_k t_k=0$,
\textit{i.e.} $P(T<+\infty)=1$. By the strong Markov's property, we can
conclude that $P(N_n-X_n>0 \, i.o.)=1$, \textit{i.e.}  $\sum_n
(N_n-X_n)=+\infty$ almost surely. Since $K_n=S_n-H_n=b+\sum_{i=1}^n
B_i(N_i-X_i)\geq \sum_{i=1}^n(N_i-X_i)$, we get $K_n=S_n-H_n\to +\infty$
almost surely.  Similarly, we can obtain that $H_n\to +\infty$ almost
surely.  \\
\indent In the general case, we have
\begin{equation*}
\begin{split}
t_k=P(T>k)&=
P(X_i=N_i\,,i=1,\dots, k)\\
&=
P(X_i=N_i\,,i=1,\dots, k\,|\, N_i\leq H_{i-1}\,, i=1,\dots,k)
P(N_i\leq H_{i-1}\,, i=1,\dots,k)\,,
\end{split}
\end{equation*}
where $P(X_i=N_i\,,i=1,\dots, k\,|\, N_i\leq H_{i-1}\,, i=1,\dots,k)$ is
equal to the product studied before and so it converges to $0$.
\end{proof}

\subsection{Almost sure convergence}
\label{subsec:as-cv}

\begin{theorem}\label{as-conv-medie-diverse}
  Assume to be in case 1) (\textit{i.e.} $m_A>m_B$). Then
  $Z_n\stackrel{a.s.}\longrightarrow Z=1$.
\end{theorem}  
\begin{proof}
  Let $e\in (m_B/m_A,\, 1)$ and set $Q_n=K_n/H_n^e$ for all $n$. Then,
  using that $(1-x)^e\leq 1-e x$ for $0\leq x\leq 1$, $H_n\leq
  H_{n+1}\leq H_n+C^2$ and~\eqref{sp-cond}, we have:
\begin{equation*}
\begin{split}
  E[Q_{n+1}/Q_n-1\,|\,\mathcal{H}_n]&=
    E\left[
      \frac{K_n+B_{n+1}(N_{n+1}-X_{n+1})}{K_n}
      \left(\frac{H_n}{H_{n+1}}\right)^e
    \,|\,\mathcal{H}_n\right]-1\\
  &=E\left[\left(\frac{H_n}{H_{n+1}}\right)^e-1\,|\,\mathcal{H}_n\right]+
  E\left[\frac{B_{n+1}(N_{n+1}-X_{n+1})}{K_n}\left(\frac{H_n}{H_{n+1}}\right)^e
    \,|\,\mathcal{H}_n\right]\\
  &\leq - e E\left[\frac{A_{n+1}X_{n+1}}{H_{n+1}}\,|\,\mathcal{H}_n\right]+
  E\left[\frac{B_{n+1}(N_{n+1}-X_{n+1})}{K_n}\,|\,\mathcal{H}_n\right]\\
  &\leq - e E\left[\frac{A_{n+1}X_{n+1}}{H_{n}+C^2}\,|\,\mathcal{H}_n\right]+
  E\left[\frac{B_{n+1}(N_{n+1}-X_{n+1})}{K_n}\,|\,\mathcal{H}_n\right]
  \\
  &= - e \frac{A_{n+1}N_{n+1}}{S_n}\frac{H_n}{H_n+C^2}+
  \frac{B_{n+1}N_{n+1}}{S_n}\,.
  \end{split}
\end{equation*}
Taking the conditional expectation with respect to $\mathcal{G}_n$ and
using~\eqref{sp-cond-rinforzi}, we get
$$
E[Q_{n+1}/Q_n-1\,|\,\mathcal{G}_n]\leq
\frac{N_{n+1}}{S_n}\left(m_{B,n+1}- e\, m_{A,n+1}\frac{H_n}{H_{n}+C^2}\right)\,.
$$ Since $H_n$ goes to $+\infty$ (see Lemma~\ref{colori-a-infinito}),
$\lim_n m_{A,n+1}=m_A>m_B=\lim_n m_{B,n+1}$ and $e\in (m_B/m_A,\,1)$,
we obtain that the above conditional expectation is smaller or equal
than zero for $n$ large enough.  It follows that, for large $n$, we
have
$$
E[Q_{n+1}-Q_n\,|\,\mathcal{G}_n]=Q_nE[Q_{n+1}/Q_n-1\,|\,\mathcal{G}_n]\leq 0
    $$ This means that $(Q_n)_n$ is eventually a positive
    (\textit{i.e.} non-negative) $\mathcal{G}$-supermartingale and so it
    converges almost surely to a finite random variable. In order to
    conclude, it is enough to observe that, since $H_n\leq S_n$,
    $S_n\stackrel{a.s.}\longrightarrow +\infty$ and $e<1$, we have
    $$
    1-Z_n=\frac{K_n}{S_n}=Q_n\frac{H_n^e}{S_n}\leq Q_n S_n^{-(1-e)}\stackrel{a.s.}
    \longrightarrow 0\,,
    $$
    that is $Z_n\stackrel{a.s.}\longrightarrow 1$.
\end{proof}

\begin{theorem}\label{lemma-base}
Assume to be in case 2). Then, we have
  \begin{equation}\label{rel-sp-cond}
|\,E[Z_{n+1}|\mathcal{G}_n]-Z_n\,|\leq
E[(A_{n+1}+B_{n+1})^2]\frac{N_{n+1}^2}{n^2}
  \end{equation}
  and so the process $(Z_n)$ is a $\mathcal{G}$-quasi-martingale and
  it almost surely converges to a random variable~$Z$ taking values in
  $[0,1]$.
\end{theorem}
It is easy to see that, in order that $(Z_n)$ is
$\mathcal{G}$-quasi-martingale, it is enough to require the condition 
\begin{equation}\label{cond-quasi-mart}
\sum_n
E[(A_{n+1}+B_{n+1})^2]\frac{E[N_{n+1}^2]}{n^2}<+\infty\,,
\end{equation}
which is obviously satisfied when $A_n\vee B_n\vee N_n\leq C$ for some
constant~$C$. Moreover, as we will see, for the proof of the above
lemma it is sufficient to assume only $m_{A,n}=m_{B,n}=m_n$ for all~$n$ (it is not necessary to have $(m_n)$ convergent).

\begin{proof}
  After some computations, we get
\begin{equation}\label{eq-base} 
Z_{n+1}-Z_n=\frac{(1-Z_n)A_{n+1}X_{n+1}-Z_nB_{n+1}(N_{n+1}-X_{n+1})}{S_{n+1}}\,.
\end{equation}

Therefore, by the model assumptions, the conditional
expectation $E[Z_{n+1}-Z_n|\mathcal{H}_n]$ is equal to
\begin{equation*}
    \sum_{k\in \mathcal{X}_{n+1}}
    [
    (1-Z_n)\frac{A_{n+1}k}{S_n+A_{n+1}k+B_{n+1}(N_{n+1}-k)}-
    Z_n\frac{B_{n+1}(N_{n+1}-k)}{S_n+A_{n+1}k+B_{n+1}(N_{n+1}-k)}
    ]
    p_{n+1,k}\,,
\end{equation*}
where $\mathcal{X}_{n+1}=\{0\vee N_{n+1}-(S_n-H_n),\dots,N_{n+1}\wedge
H_n\}$ and $p_{n+1,k}=p_k(N_{n+1},S_n,H_n)$ is given by
\eqref{hypergeometric-distr}. We observe that $\mathcal{X}_{n+1}$ and
$p_{n+1,k}$ are $\mathcal{G}_n$-measurable and so the conditional
expectation $E[Z_{n+1}-Z_n|\mathcal{G}_n]$ is equal to
\begin{equation*}
    \sum_{k\in \mathcal{X}_{n+1}} \left\{ (1-Z_n)
    E\left[\frac{A_{n+1}k}{S_n+A_{n+1}k+B_{n+1}(N_{n+1}-k)}|\mathcal{G}_n\right]-
    Z_n
    E\left[\frac{B_{n+1}(N_{n+1}-k)}{S_n+A_{n+1}k+B_{n+1}(N_{n+1}-k)}|\mathcal{G}_n
      \right]
    \right\}p_{n+1,k}\,.
\end{equation*}
Now, we consider the above quantity and we add and subtract the
quantity $A_{n+1}k/S_n$ in the first conditional expectation and the
quantity $B_{n+1}(N_{n+1}-k)/S_n$ in the second conditional
expectation, so that the two conditional expectations can be rewritten
respectively as
\begin{equation*}
  \begin{split}
    &E\left[
  \frac{-A_{n+1}^2k^2-A_{n+1}B_{n+1}k(N_{n+1}-k)}{S_n[S_n+A_{n+1}k+B_{n+1}(N_{n+1}-k)]}
  |\mathcal{G}_n\right]+\frac{m_nk}{S_n}\\
    &E\left[
      \frac{-B_{n+1}^2(N_{n+1}-k)^2-A_{n+1}B_{n+1}k(N_{n+1}-k)}
           {S_n[S_n+A_{n+1}k+B_{n+1}(N_{n+1}-k)]}
  |\mathcal{G}_n\right]+\frac{m_n(N_{n+1}-k)}{S_n}\,,
  \end{split}
\end{equation*}
where we have used \eqref{sp-cond-rinforzi} and the fact that
$m_{A,n}=m_{B,n}=m_n$. Finally, we observe that
$$
\sum_{k\in \mathcal{X}_{n+1}} \frac{(1-Z_n)m_n k - Z_n m_n(N_{n+1}-k)}{S_n} p_{n+1,k}=
\frac{m_n}{S_n}\sum_{k\in \mathcal{X}_{n+1}}(k-N_{n+1}Z_n)p_{n+1,k}=0,
$$ because $\sum_{k\in \mathcal{X}_{n+1}}kp_{n+1,k}$ is the mean value of
the hypergeometric distribution with parameters $N_{n+1},\, S_n,\,
H_n$ and so it is equal to $N_{n+1}H_n/S_n=N_{n+1}Z_n$.  Summing up,
the conditional expectation $E[Z_{n+1}-Z_n|\mathcal{G}_n]$ is equal to
$$
\sum_{k\in \mathcal{X}_{n+1}}
E\left[
  \frac{Z_nB_{n+1}^2(N_{n+1}-k)^2-(1-Z_n)A_{n+1}^2k^2
    +(2Z_n-1)A_{n+1}B_{n+1}k(N_{n+1}-k)}
       {S_n[S_n+A_{n+1}k+B_{n+1}(N_{n+1}-k)]}
  |\mathcal{G}_n\right]p_{n+1,k}\,.
$$
Therefore, using assumption (A3), we have 
$$
|\,E[Z_{n+1}|\mathcal{G}_n]-Z_n\,|\leq
E\left[\frac{(A_{n+1}+B_{n+1})^2N_{n+1}^2}{S_n^2}|\mathcal{G}_n\right]
=E[(A_{n+1}+B_{n+1})^2]\frac{N_{n+1}^2}{S_n^2}
$$ and, since $A_n\wedge B_n\wedge N_n\geq 1$ by definition, we
finally get \eqref{rel-sp-cond}. When condition
\eqref{cond-quasi-mart} is satisfied (as when $A_n\vee B_n\vee N_n\leq
C$ for some constant $C$), the process $(Z_n)$ is a
$\mathcal{G}$-martingale taking values in $[0,1]$ and, hence, it
almost surely converges to some random variable $Z$ taking values in
$[0,1]$.
\end{proof}

\begin{remark}
\rm From \eqref{eq-base}, we immediately get that, if $A_n=B_n$ for
all $n$, then
\begin{equation*}
  Z_{n+1}-Z_n=\frac{A_{n+1}(X_{n+1}-Z_nN_{n+1})}{S_{n}+A_{n+1}N_{n+1}}
\end{equation*}
and so $(Z_n)$ is an $\mathcal{H}$-martingale, because of assumptions
(A1) and (A2). Therefore, for its almost sure convergence, it is not
necessary condition \eqref{cond-quasi-mart}.  This is the case
considered in~\cite{perron,crimaldi2016}.
\end{remark}

\begin{remark}\rm 
Lemma~\ref{lemma-app} (with $Y_n=X_n/N_n$) immediately
implies that, in both cases 1) and 2),  the sequence
\begin{equation}\label{def-M}
M_n=\frac{1}{n}\sum_{j=1}^{n}\frac{X_j}{N_j},
\end{equation}
which is the empirical mean of the proportion, in the samples, of
balls of color~A, also converges almost surely to~$Z$.
\end{remark}

\begin{proposition}\label{conv-S}
  Assume to be in one of the previous two cases 1) and 2) and let
  $Z\stackrel{a.s.}=\lim_n Z_n$. Moreover, assume
  \begin{equation}\label{ass-base-2}
  E[N_n|{\mathcal F}_{n-1}]\stackrel{a.s.}\longrightarrow N\,,
  \end{equation}
  where $N$ is a (strictly positive finite) random variable.\\
  \indent Then
$$
\frac{H_n}{n}\stackrel{a.s.}\longrightarrow m_A N Z\,,\qquad
\frac{K_n}{n}=\frac{S_n-H_n}{n}\stackrel{a.s.}\longrightarrow m_B N (1-Z).
$$
and so   
$$
\frac{S_n}{n}\stackrel{a.s.}\longrightarrow m_A N Z + m_B N (1-Z).
$$
\end{proposition}
%
%
\begin{proof} It is enough to apply Lemma~\ref{lemma-app}
  with $Y_j=A_jX_j$ (resp. $Y_j=B_j(N_{j}-X_j)$. Indeed,
  we have $Y_j\leq A_jN_j$ (resp. $Y_j\leq B_jN_j$) for each $j$ and so 
  $E[Y_j^2]\leq E[(A_j+B_j)^2]E[N_j^2]$.  Moreover  
  \begin{equation*}
    \begin{split}
  E[A_jX_j |{\mathcal F}_{j-1}]&=
E\left[
    E[\,
      E[A_jX_j|\mathcal{H}_{j-1}]\,
      |\mathcal{G}_{j-1}]\,
    |\mathcal{F}_{j-1}\right]
  =E\left[\,
    E[A_jN_jZ_{j-1}|\mathcal{G}_{j-1}]
    |{\mathcal F}_{j-1}\right]\\
  &=E[m_{A,j}N_jZ_{j-1}|{\mathcal F}_{j-1}]
= m_{A,j}E[N_j|\mathcal{F}_{j-1}]Z_{j-1}\stackrel{a.s.}\longrightarrow m_A N Z
    \end{split}
  \end{equation*}
and 
\begin{equation*}
  \begin{split}
E[B_j(N_j-X_j) |{\mathcal F}_{j-1}]&=
  E\left[
    E[\,
      E[B_j(N_j-X_j)|\mathcal{H}_{j-1}]\,
      |\mathcal{G}_{j-1}]\,
    |\mathcal{F}_{j-1}\right]\\
  &=E\left[\,
    E[B_jN_j(1-Z_{j-1})|\mathcal{G}_{j-1}]
    |{\mathcal F}_{j-1}\right]\\
  &=E[m_{B,j}N_j(1-Z_{j-1})|{\mathcal F}_{j-1}]\\
  &=m_{B,j}E[N_j|\mathcal{F}_{j-1}](1-Z_{j-1})
\stackrel{a.s.}\longrightarrow m_B N(1-Z)\,.
    \end{split}
  \end{equation*}
Therefore, we have $H_n/n\stackrel{a.s.}\longrightarrow m_ANZ$ and
$K_n/n\stackrel{a.s.}\longrightarrow m_B N(1-Z)$ and so
$S_n/n=H_n/n+K_n/n\stackrel{a.s.}\longrightarrow m_ANZ+m_BN(1-Z)$.
\end{proof}

\begin{remark}\label{rem-caso-1}
  When we are in case 1), then $Z=1$ almost surely and so we have
  $H_n$ and $S_n$ go to $+\infty$ with rate $n$. Moreover, we observe
  that, for each $e\in (m_B/m_A,\,1)$, we have
  $$
  n^{1-e}(1-Z_n)=n^{1-e}\frac{K_n}{S_n}=
  \left(\frac{n}{S_n}\right)^{1-e}\left(\frac{H_n}{S_n}\right)^eQ_n\,,
  $$ where $Q_n$ is defined as in the proof of Theorem
 ~\ref{as-conv-medie-diverse}. Since $n/S_n$, $H_n/S_n$ and $Q_n$
  converge almost surely to suitable finite random variables, we get
  that $n^{1-e}(1-Z_n)$ converges almost surely to a finite random
  variable. Since $e$ is arbitrary, we necessarily have
  $n^{1-e}(1-Z_n)\stackrel{a.s.}\longrightarrow 0$, that is, for all
  $e\in (m_B/m_A,\,1)$, we have $1-Z_n\stackrel{a.s.}{=}o(n^{-(1-e)})$
  and so $K_n=S_n(1-Z_n)=o(n^e)$. 
 \end{remark}

When we are in case 2), since $m N >0$ almost surely, the above limit
result implies that $S_n$ goes to $+\infty$ with rate $n$; while it is
not sufficient in order to get some information on the asymptotic
behavior of $H_n$ and $K_n$, because~$Z$ may assume the value $0$ or
$1$.  In the sequel, we will prove that both $H_n$ and $K_n$ go to
$+\infty$ at rate $n$.

\begin{theorem}\label{th-legge-1}
  Assume to be in case 2) and assume condition
  \eqref{ass-base-2}. Then we have $P(Z=0)+P(Z=1)=0$. (Consequently
  the rate at which $H_n$ and $K_n$ go to $+\infty$ is equal to~$n$.)
\end{theorem}
\begin{proof} 
  Set $Y_n=\ln(H_n/K_n)$, $\Delta_n=E[Y_{n+1}-Y_n|\mathcal{G}_n]$ and
  $Q_n=E[(Y_{n+1}-Y_n)^2]$. If we prove $\sum_n\Delta_n<+\infty$ and
  $\sum_n Q_n<+\infty$ almost surely, then $Y_n$ converges almost
  surely to a finite random variable (see Lemma~3.2
  in~\cite{pemantle-volkov-1999}). This fact implies that $H_n/K_n$
  converges to a random variable $Y$ with values in $(0,+\infty)$. It
  follows that $Z_n=\frac{H_n}{S_n}=\frac{H_n/K_n}{H_n/K_n+1}$
  converges almost surely to $Y/(Y+1)$, which is a random variable
  with values in $(0,1)$. Then $P(Z=0)+P(Z=1)=0$.  \\ \indent The rest
  of the proof is devoted to verify that $\sum_n\Delta_n<+\infty$ and
  $\sum_n Q_n<+\infty$ almost surely.\\ To this regard, we recall
  that, by Lemma~\ref{rate-colori-2}, we have $1/K_n=O(1/n^\gamma)$
  and $1/H_n=O(1/n^\gamma)$ with $\gamma>0$. Moreover, using the
  notation \eqref{hypergeometric-distr}, we have
  \begin{equation*}
  \begin{split}
    &E[\ln(H_{n+1})-\ln(H_n)|\mathcal{H}_n]-E[\ln(K_{n+1})-\ln(K_n)|\mathcal{H}_n]=
    \\
    &\sum_{k\in\mathcal{X}_{n+1}}
    \left\{\left(\ln(H_{n}+A_{n+1}k)-\ln(H_n)\right)-
    \left(\ln(K_{n}+B_{n+1}(N_{n+1}-k))-\ln(K_n)\right)\right\}p_{n+1,k}=
    \\
   &\sum_{k\in\mathcal{X}_{n+1}}
    \left\{\int_0^{A_{n+1}k}\frac{1}{H_n+t}\,dt-
    \int_0^{B_{n+1}(N_{n+1}-k)}\frac{1}{K_n+t}\,dt\right\}p_{n+1,k} 
  \end{split}
  \end{equation*}
  Since $1/(H_n+t)\leq 1/H_n$ and $1/(K_n+t)\geq 1/K_n-t/K_n^2$ for
  each $t\geq 0$ and each $n$, the last term of the above equalities
  is eventually smaller or equal than
  $$
\sum_{k\in\mathcal{X}_{n+1}}\left\{\frac{A_{n+1}k}{H_n}-
\frac{B_{n+1}(N_{n+1}-k)}{K_n}+
c\frac{B_{n+1}^2(N_{n+1}-k)^2}{2 K_n^2}\right\}p_{n+1,k}\,.
$$
Now, we observe that
$$
E[\sum_{k\in\mathcal{X}_{n+1}}\left(\frac{A_{n+1}k}{H_n}-
  \frac{B_{n+1}(N_{n+1}-k)}{K_n}\right)p_{n+1,k}\,|\,\mathcal{G}_n]=
\frac{m_{n+1}N_{n+1}}{S_n}-\frac{m_{n+1}N_{n+1}}{S_n}=0\,.
$$ Therefore, we have for $n$ large enough (using $(1-Z_n)=K_n/S_n$)
$$
\Delta_n\leq \frac{cC^2}{2 K_n^2}
\left\{Z_n(1-Z_n)N_{n+1}\frac{S_n-N_{n+1}}{S_n-1}+N_{n+1}^2(1-Z_n)^2\right\}
=O(1/(K_nS_n))=O(1/n^{1+\gamma})\,.
$$
Similarly, we have
 \begin{equation*}
  \begin{split}
&E[(\ln(H_{n+1})-\ln(H_n)-\ln(K_{n+1})+\ln(K_n))^2|\mathcal{H}_n]
\leq\\
&2\left\{
E[(\ln(H_{n+1})-\ln(H_n))^2|\mathcal{H}_n]+
E[(\ln(K_{n+1})-\ln(K_n))^2|\mathcal{H}_n]
\right\}
\leq\\
&2
\sum_{k\in\mathcal{X}_{n+1}}
\left(\frac{A_{n+1}^2k^2}{H_n^2}+\frac{B_{n+1}^2(N_{n+1}-k)^2}{K_n^2}\right)p_{n+1,k}
=O(1/(H_nS_n))+O(1/K_nS_n)=O(1/n^{1+\gamma})\,.
  \end{split}
 \end{equation*}
 The last statement (into the brackets) immediately follows from Proposition 
~\ref{conv-S}.
\end{proof}


\subsection{Central limit theorems for the case of equal
  reinforcement means}
\label{subsec:CLT}

Since in case~2), the limit proportion is a random variable~$Z$, in
the sequel we provide results in order to get some information on it.

\begin{theorem}\label{main-1}
Assume to be in case~2) and assume condition
\eqref{ass-base-2}. Moreover, suppose to have 
\begin{equation}\label{ass-N-sp-cond-2}
E[N_n^2|\mathcal{F}_{n-1}]\stackrel{a.s.}\longrightarrow Q\,,
\end{equation}
where $Q$ is a (strictly positive finite) random variable, 
and
\begin{equation}\label{ass-A-B}
  q_{A,n}=E[A_n^2]\to q_A\,,\quad q_{B,n}=E[B_n^2]\to q_B\,,
  \quad q_{AB,n}=E[A_nB_n]\to q_{AB}\,,
\end{equation}
where $q_A,\,q_B$ and $q_{AB}$ are (strictly positive finite) constants. 
\\ \indent Then
$\sqrt{n}(Z_n-Z)$ converges in the sense of the almost sure
conditional convergence with respect to ${\mathcal F}=({\mathcal F}_n)$
to the Gaussian kernel ${\mathcal N}(0, V)$, where
\begin{equation}\label{eq-V}
  \begin{split}
V&=Z(1-Z)
\frac{(1-Z)q_A[(1-Z)N+ZQ]+Zq_B[ZN+(1-Z)Q]-2Z(1-Z)q_{AB}(Q-N)}
     {(mN)^2}\\
     &=Z(1-Z)
     \frac{N[(1-Z)^2q_A+Z^2q_B+2Z(1-Z)q_{AB}]+Z(1-Z)Q[q_A+q_B-2q_{AB}]}
          {(mN)^2}\,.
          \end{split}
\end{equation}
\end{theorem}
Before the proof, we premise some remarks.

\begin{remark}\label{rem-on-V} \rm 
Regarding formula \eqref{eq-V}, recall that $N\geq 1$ a.s., $Q\geq 1$
a.s., $q_A\geq 1\,q_B\geq 1,\,q_{AB}\geq 1$ and
$q_A+q_B-2q_{AB}=\lim_n E[(A_n-B_n)^2]\geq 0$. Moreover, we have
proven that $P(Z=0)=P(Z=1)=0$ (see Theorem
\ref{th-legge-1}). Therefore, we have $P(V>0)=1$. In addition, we note
that $V$ is not degenerate provided $P(Z=z)<1$ for all $z\in
(0,1)$. For this last fact, we refer to the next
Theorem~\ref{th-legge-2}, which states that we also have $P(Z=z)=0$
for all $z\in (0,1)$.
\end{remark}

\begin{remark} \rm When $A_n=B_n$ for all $n$, we have $q_A=q_B=q_{AB}=q$ and so 
   we get $V=Z(1-Z)q/(m^2 N)$, that does not depend on $Q$. Indeed, in
   this case the above assumption \eqref{ass-N-sp-cond-2} can be
   deleted (see~\cite{crimaldi2016}).
\end{remark}

\begin{remark} \rm When $N_n=k$ for each $n$, with $k$ a fixed constant,
  we have
  \begin{equation}\label{eq-V-caso-costante}
    \begin{split}
V&=kZ(1-Z)
\frac{(1-Z)^2q_A+Z^2q_B+2Z(1-Z)q_{AB}+Z(1-Z)k(q_A+q_B-2q_{AB})}
     {(mk)^2}\\
    & =Z(1-Z)
\frac{(1-Z)^2q_A+Z^2q_B+Z(1-Z)[k(q_A+q_B)-2q_{AB}(k-1)]}
     {m^2k}\,.
     \end{split}
\end{equation}
  In particular, for $k=1$,
%
%
 we observe that $V$ does not depend on $q_{AB}$.
\end{remark}

\begin{proof} Setting $X'_n=X_n/N_n$
for each $n$, the sequence $(X'_n)$ is $\mathcal G$-adapted and
bounded. Moreover, we have
\begin{equation}\label{eq-p0}
\begin{split}
E[X_{n+1}'|{\mathcal G}_n]=
E[ N_{n+1}^{-1} X_{n+1}|{\mathcal G}_n]=  N_{n+1}^{-1} E[X_{n+1}|{\mathcal G}_n]=
N_{n+1}^{-1} N_{n+1} Z_n=Z_n \,.
\end{split}
\end{equation}
We want to apply Theorem~\ref{th-main-app} to $Y_n=X_n'$. By Theorem 
~\ref{lemma-base}, we have
$$
n^3 E\left[\,(E[Z_{n+1}|\mathcal{G}_n]-Z_n)^2\,\right]\longrightarrow 0.
$$ Therefore, in order to prove Theorem~\ref{main-1}, it suffices to
prove that the following conditions are satisfied
\begin{itemize}
\item[c1)] $E[\sup_{j\geq 1} \sqrt{j} |Z_{j-1}-Z_j|\,]<+\infty$;
\item[c2)] $n\sum_{j\geq n}(Z_{j-1}-Z_j)^2\stackrel{a.s.}\longrightarrow V$.
\end{itemize}
In the following we verify the above conditions.\\
\noindent{\em Condition c1).} We observe that, by \eqref{eq-base} and
recalling that $A_j\wedge B_j\wedge N_j\geq 1$ and $A_j\vee B_j\vee
N_j\leq C$, we have
\begin{equation}\label{eq-p2}
|Z_{j-1}-Z_j|\leq
\frac{(A_{j}+B_{j})N_{j}}{j}\leq \frac{2C^2}{j}\,.
\end{equation}
Therefore condition c1) is obviously verified.
\\
\noindent{\em Condition c2).} We want to apply Lemma~\ref{lemma-app}
with $Y_j=j^2(Z_{j-1}-Z_j)^2$.  By the assumptions and
inequality~\eqref{eq-p2}, we have $\sum_{j\geq 1}
j^{-2}E[Y_j^2]<+\infty$. Moreover, by equality~\eqref{eq-base}, we
have
\begin{equation*}
  (Z_{j-1}-Z_j)^2=
  \frac{(1-Z_{j-1})^2A_j^2N_j^2(X_j')^2}{S_j^2}
  +
  \frac{Z_{j-1}^2B_j^2N_j^2(1-X_j')^2}{S_j^2}
  -2
  \frac{Z_{j-1}(1-Z_{j-1})A_jB_jN_j^2X_j'(1-X_j')}{S_j^2}\,.
\end{equation*}
Therefore, we study the convergence of the following three terms:
\begin{itemize}
\item 
  $T_{1,j-1}=j^2
  E\left[\frac{(1-Z_{j-1})^2A_j^2N_j^2(X_j')^2}{S_j^2}|\mathcal{F}_{j-1}\right]$,
\item $T_{2,j-1}=j^2
  E\left[\frac{Z_{j-1}^2B_j^2N_j^2(1-X_j')^2}{S_j^2}|\mathcal{F}_{j-1}\right]
  $,
  \item
    $T_{3,j-1}=j^2E\left[\frac{Z_{j-1}(1-Z_{j-1})A_jB_jN_j^2X_j'(1-X_j')}{S_j^2}|
    \mathcal{F}_{j-1}\right]$.
\end{itemize}
Consider the first term $T_{1,j-1}$. By assumption (A3), we get the
two inequalities:
\begin{equation*}
\begin{split}
& T_{1,j-1}\geq 
\frac{j^2}{(S_{j-1} + C^2)^2} (1-Z_{j-1})^2
E[A_j^2] E[N_j^2 (X'_j)^2|{\mathcal F}_{j-1}]\\
& 
T_{1,j-1} \leq 
\frac{j^2}{S_{j-1}^2} (1-Z_{j-1})^2E[A_j^2] E[N_j^2 (X'_j)^2|{\mathcal F}_{j-1}].
\end{split}
\end{equation*}
Since $S_n/n\stackrel{a.s.}\longrightarrow N m>0$,
$Z_{j-1}\stackrel{a.s.}\longrightarrow Z$ and $E[A_j^2]\to q_A$, it is
enough to verify the almost sure convergence of $E[N_j^2
  (X'_j)^2|{\mathcal F}_{j-1}]$. To this purpose, we
observe that we can write
\begin{equation*}
E[N_j^2 (X'_j)^2|{\mathcal F}_{j-1}]=
E\left[ N_j^2 
E[(X'_j)^2 |{\mathcal G}_{j-1}]
\,|\,{\mathcal F}_{j-1}\right]
\end{equation*}
and, by (A2), the conditional expectation
$E[(X'_j)^2 |{\mathcal G}_{j-1}]$ coin\-ci\-des with
\begin{equation*}
\begin{split}
N_j^{-2} E[X_j^2|{\mathcal G}_{j-1}]
&=N_j^{-2}
\left[
Z_{j-1}(1-Z_{j-1})(S_{j-1}-1)^{-1} N_{j}(S_{j-1}-N_j)
+Z_{j-1}^2 N_j^2\right]\\
&=
Z_{j-1}(1-Z_{j-1})(S_{j-1}-1)^{-1}N_j^{-1}\left(S_{j-1}-N_j\right)
+Z_{j-1}^2.
\end{split}
\end{equation*}
Therefore we obtain 
\begin{equation*}
E[N_j^2 (X'_j)^2|{\mathcal F}_{j-1}]
=Z_{j-1}(1-Z_{j-1})(S_{j-1}-1)^{-1}
\left(S_{j-1} E[N_j|{\mathcal F}_{j-1}]-E[N_j^2|\mathcal{F}_{j-1}]\right)
+Z_{j-1}^2E[N_j^2|\mathcal{F}_{j-1}],
\end{equation*}
which converges almost surely to $Z(1-Z)N+Z^2Q$ (since $E[N_j^2|{\mathcal
    F}_{j-1}]$ is bounded by $C^2$ and
$S_{j-1}\stackrel{a.s.}\longrightarrow +\infty$).  Hence
$T_{1,j-1}$ converges almost surely to
$T_1=Z(1-Z)^2q_A(mN)^{-2}[(1-Z)N+ZQ]$.
Similarly, we get
\begin{equation*}
  \begin{split}
E[N_j^2 (1-X'_j)^2|{\mathcal F}_{j-1}]&=
E[N_j^2|\mathcal{F}_{j-1}]+E[N_j^2(X_j')^2|\mathcal{F}_{j-1}]
-2E[N_j^2X_j'|\mathcal{F}_{j-1}]\\
&=E[N_j^2|\mathcal{F}_{j-1}]+E[N_j^2(X_j')^2|\mathcal{F}_{j-1}]
-2Z_jE[N_j^2|\mathcal{F}_{j-1}]
\\
&\longrightarrow Q+Z(1-Z)N+Z^2Q-2ZQ=Z(1-Z)N+(1-Z)^2Q.
  \end{split}
  \end{equation*}
and so 
$T_{2,j-1}$ converges almost surely to
$T_2=Z^2(1-Z)q_B(mN)^{-2}[ZN+(1-Z)Q]$. Finally, we have
\begin{equation*}
  \begin{split}
E[N_j^2 X_j'(1-X'_j)|{\mathcal F}_{j-1}]&=
E[N_j^2X_j'|\mathcal{F}_{j-1}]-E[N_j^2(X_j')^2|\mathcal{F}_{j-1}]
\\
&=Z_{j-1}E[N_j^2|\mathcal{F}_{j-1}]-E[N_j^2(X_j')^2|\mathcal{F}_{j-1}]
\\
&\longrightarrow ZQ-Z(1-Z)N-Z^2Q=Z(1-Z)(Q-N).
  \end{split}
  \end{equation*}
and so $T_{3,j-1}$ converges almost surely to
$T_3=Z^2(1-Z)^2q_{AB}(mN)^{-2}(Q-N)$. By Lemma~\ref{lemma-app}, condition c2)
is satisfied with $V=T_1+T_2-2T_3$.  The proof is so
concluded.
\end{proof}

\begin{theorem}\label{main-2}
Under the assumptions of Theorem~\ref{main-1}, suppose also that  
\begin{equation}\label{ass-base-bis}
E[N_n^{-1}|{\mathcal F}_{n-1}]\stackrel{a.s.}\longrightarrow L\,,
\end{equation}
where $L$ is a (positive bounded) random variable.\\ 
\indent Then 
$$
[ \sqrt{n}(M_n-Z_n), \sqrt{n}(Z_n-Z) ] 
\stackrel{stably}\longrightarrow 
{\mathcal N}(0, U)\otimes {\mathcal N}(0, V),
$$
where $M_n$ is defined in \eqref{def-M}, $V$ is defined in
\eqref{eq-V} and $U=V+Z(1-Z)[L-2N^{-1}]$.
\end{theorem}
In particular, we have that $\sqrt{n}(M_n-Z_n)$ converges stably to
${\mathcal N}(0, U)$ and $\sqrt{n}(M_n-Z)$ converges stably to
${\mathcal N}(0, U+V)$, with $U+V>0$ almost surely (see Remark
\ref{rem-on-V}).

\begin{remark}\label{rem-on-U}
  \rm Regarding the limit random variance $U$, we note that, by Jensen
  inequality, we have $(E[N_n|{\mathcal F}_{n-1}])^2\leq
  E[N_n^2|\mathcal{F}_{n-1}]$ and $E[N_n|\mathcal{F}_{n-1}]^{-1}\leq
  E[N_n^{-1}|\mathcal{F}_{n-1}]$ and so we have $N^2\leq Q$ and $1/N\leq L$.
  Therefore, we get
  \begin{equation*}
    \begin{split}
 & V\geq Z(1-Z)
     \frac{N[(1-Z)^2q_A+Z^2q_B+2Z(1-Z)q_{AB}]+Z(1-Z)N^2[q_A+q_B-2q_{AB}]}
          {(mN)^2}
          \quad\mbox{and}\\
          \quad &L-\frac{2}{N}\geq -\frac{1}{N}\,.
          \end{split}
\end{equation*}
          Moreover, since $N_n\geq 1$ for each $n$, we have $N\geq 1$
          and so $N^2\geq N$. It follows the relation $V\geq
          Z(1-Z)[(1-Z)q_A+Zq_B]/(mN)^2$ and hence
          $$
          U\geq \frac{Z(1-Z)}{N}\left[\frac{(1-Z)q_A+Zq_B}{m^2}-1\right]\,.
          $$ Since $q_A\geq m^2$ and $q_B\geq m^2$ and
          $P(Z=0)=P(Z=1)=0$, the quantity in the right side of the
          last inequality is always greater or equal than zero almost
          surely and it is equal to zero if and only if
          $q_A=q_B=m^2$. Summing up, the rate of convergence of
          $(M_n-Z_n)$ to zero is $1/2$ whenever $q_A>m^2$ or $q_B>m^2$
          and, otherwise, it could be even greater.
\end{remark}

\begin{proof} Thanks to what we have already
proven in the previous proof, it suffices to verify that the following
condition is satisfied (see Theorem~\ref{th-main-app}
applied to $Y_n=X_n'$):
\begin{itemize}
\item[c3)] $n^{-1} \sum_{j=1}^n
  \big[X'_j-Z_{j-1}+j(Z_{j-1}-Z_j)\big]^2
  \stackrel{P}\longrightarrow U$.
\end{itemize}
To this purpose, we apply Lemma~\ref{lemma-app} with
$$Y_j=\big[ X'_j - Z_{j-1}+ j(Z_{j-1}-Z_j) \big]^2.$$ Indeed, by the
assumptions, we have $\sum_{j\geq 1}
j^{-2}E[Y_j^2]<+\infty$. Moreover, from what we have already seen in
the previous proof, we can get
\begin{equation*}
j^2E[(Z_{j-1} - Z_{j})^2|{\mathcal F}_{j-1}]
\stackrel{a.s.}\longrightarrow V\,.
\end{equation*}
Moreover, leveraging the above computations, we have
\begin{equation*}
  \begin{split}
  E[(X'_j-Z_{j-1})^2|{\mathcal F}_{j-1}]
&= E[(X'_j)^2|{\mathcal F}_{j-1}]-Z_{j-1}^2
  \\
  &= Z_{j-1}(1-Z_{j-1})(S_{j-1}-1)^{-1}
  \left(S_{j-1}E[N_j^{-1}|\mathcal{F}_{j-1}]-1\right)
  \stackrel{a.s.}\longrightarrow
  Z(1-Z)L\,.
  \end{split}
\end{equation*}
Finally, we observe that 
\begin{equation*}
\begin{split}
&j(X_j'-Z_{j-1})(Z_{j-1}-Z_j)=-
j(X_j'-Z_{j-1})\frac{(1-Z_{j-1})A_{j}N_jX'_{j}-Z_{j-1}B_{j}N_j(1-X'_{j})}{S_{j}}
=\\
&  -\frac{j(1-Z_{j-1})A_{j}N_j(X'_{j})^2}{S_{j}}
    +\frac{jZ_{j-1}(1-Z_{j-1})A_{j}N_jX'_{j}}{S_j}
    +\frac{jZ_{j-1}B_{j}N_jX_j'(1-X'_{j})}{S_j}
    -\frac{jZ_{j-1}^2B_{j}N_j(1-X'_{j})}{S_{j}}
    =\\
 & -U_{1,j}+U_{2,j}+U_{3,j}-U_{4,j}\,.
\end{split}
\end{equation*}
With the same techniques adopted in the previous proof, we can get
\begin{equation*}
  \begin{split}
  &T_{1,j-1}=E[U_{1,j}|\mathcal{F}_{j-1}]\stackrel{a.s.}\longrightarrow
 T_1=Z(1-Z)^2/N+Z^2(1-Z)
  \\
  &T_{2,j-1}=E[U_{2,j}|\mathcal{F}_{j-1}]\stackrel{a.s.}\longrightarrow
  T_2=Z^2(1-Z)
  \\
  &T_{3,j-1}=E[U_{3,j}|\mathcal{F}_{j-1}]\stackrel{a.s.}\longrightarrow
  T_3=Z^2-Z^2(1-Z)/N-Z^3=-Z^2(1-Z)/N+Z^2(1-Z)
  \\
  &T_{4,j-1}=E[U_{4,j}|\mathcal{F}_{j-1}]\stackrel{a.s.}\longrightarrow
  T_4=Z^2(1-Z)
  \end{split}
\end{equation*}
Summing up, we obtain the almost sure convergence of
$E[Y_j|\mathcal{F}_{j-1}]$ to $U=V+Z(1-Z)L+2(-T_1+T_2+T_3-T_4)=
V+Z(1-Z)(L-2N^{-1})$.
\end{proof}

\subsection{Probability distribution of the limit proportion
  in the case of equal reinforcement means}
\label{subsec:no-atom}

When we are in case 2), the distribution of the limit proportion $Z$
is unknown except in a few particular cases (see~\cite{perron}).  What
we are able to prove in the general case is that it is diffuse (see
Theorem~\ref{th-legge-2} below) and to leverage the above central
limit theorems in order to get asymptotic confidence intervals for~$Z$
(see Subsection~\ref{conf-int} below).

\begin{theorem}\label{th-legge-2} 
  Assume the same assumptions as in Theorem~\ref{main-1}, then
  $P(Z=z)=0$ for all $z\in [0,1]$.
\end{theorem}

\begin{proof}
We already know that $P(Z=0)=P(Z=1)=0$ (see Theorem~\ref{th-legge-1})
In order to prove that $P(Z=z)=0$ for all $z\in (0,1)$, we can argue
exactly as done in~\cite[Cor. 4.1]{crimaldi2016} or in Th.~3.2 
in~\cite{cri-dai-min}.  Since
the key issue on which the proof is based is the almost sure
conditional convergence of $\sqrt{n}(Z_n-Z)$ with respect to
${\mathcal F}=({\mathcal F}_n)$ to a Gaussian kernel ${\mathcal N}(0,
V)$, for some $V>0$ on $\{Z\in (0,1)\}$.
\end{proof}

\subsection{Asymptotic confidence intervals for the limit proportion in the case of equal reinforcement means} 
\label{conf-int}

Suppose to be in case 2). By means of Theorem~\ref{main-1} and Theorem
~\ref{main-2} (together with Theorem \ref{lemma-conv-stab}), we can
construct \emph{asymptotic confidence intervals} for the limit
proportion $Z$. More precisely, assume $A_n\vee B_n\vee N_n\leq C$ for
each $n$ and conditions~\eqref{ass-base-2}, \eqref{ass-N-sp-cond-2},
and~\eqref{ass-A-B}. Then, by Lemma \ref{lemma-app}, the random variables 
\begin{equation}\label{eq-estimators}
  \widehat{m}_{n}=\frac{\sum_{j=1}^n A_j}{n},\qquad
  \widehat{q}_{A,n}=\frac{\sum_{j=1}^n A_j^2}{n}, \qquad
  \widehat{q}_{B,n}=\frac{\sum_{j=1}^n B_j^2}{n},
  \qquad \widehat{q}_{AB,n}=\frac{\sum_{j=1}^n A_jB_j}{n}
\end{equation}
are strongly consistent estimators of the constants $m,\,q_A,\,q_B$
and $q_{AB}$ (supposed unknown), respectively. By Lemma
\ref{lemma-app} again, the random variables
\begin{equation}
  \widehat{\mu}_n=\frac{\sum_{j=1}^n N_j}{n},\qquad
  \widehat{q}_{N,n}=\frac{\sum_{j=1}^n N^2_j}{n}\,,
\end{equation}
are strongly consistent estimators of the random variables $N$ and $Q$. 
Hence, the random variable
\begin{equation*}
  \begin{split}
  &V_n=Z_n(1-Z_n)\times \\ &
    \frac{(1-Z_n)\widehat{q}_{A,n}[(1-Z_n)\widehat{\mu}_n+Z_n\widehat{q}_{N,n}]+
      Z_n\widehat{q}_{B,n}[Z_n\widehat{\mu}_n+(1-Z_n)\widehat{q}_{N,n}]
      -2Z_n(1-Z_n)\widehat{q}_{AB,n}(\widehat{q}_{N,n}-\widehat{\mu}_n)}
         {(\widehat{m}_n\widehat{\mu}_n)^2}
       \end{split}
\end{equation*}
results a strongly consistent estimator of the random variable $V$
(defined in Theorem~\ref{main-1}).  Recalling that $V>0$ almost surely
(see Remark \ref{rem-on-V}), by Theorem~\ref{main-1}, together with
Theorem \ref{lemma-conv-stab}, we obtain that a confidence interval
for $Z$ is
\begin{equation} \label{conf-int-formulae-1}
Z_n\pm q_{1-\frac{\alpha}{2}}\sqrt{\frac{V_n}{n}}\,,
\end{equation}
where $q_{1-\frac{\alpha}{2}}$ is the quantile of order
$1-\frac{\alpha}{2}$ of the standard normal distribution.\\ \indent
When $N_n=k$, with $k$ a known constant, for $V_n$ we can employ the
simpler formula \eqref{eq-V-caso-costante} with
$\widehat{q}_{A,n},\,\widehat{q}_{B,n}$ and $\widehat{q}_{AB,n}$
instead of $q_A,\,q_B$ and $q_{AB}$.  \\ \indent If
condition~\eqref{ass-base-bis} is also satisfied, then, again by Lemma
\ref{lemma-app}, $\widehat{\eta}_n=\frac{\sum_{j=1}^n N_j^{-1}}{n}$ is
a strongly consistent estimator of the random variable $L$ (defined in
Theorem~\ref{main-2}) and so, setting
\begin{equation*}
W_n=2V'_n+M_n(1-M_n)[\widehat{\eta}_n-2(\widehat{\mu}_n)^{-1}]\,,
\end{equation*}
where $V_n'$ is equal to $V_n$ but with $M_n$ instead of $Z_n$, is a
strongly consistent estimator of the random variable $W=U+V$.  Since
$U+V>0$ almost surely, by Theorem~\ref{main-2}, together with Theorem
\ref{lemma-conv-stab}), we get that
\begin{equation} \label{conf-int-formulae-2}
M_n\pm q_{1-\frac{\alpha}{2}}\sqrt{\frac{W_n}{n}}
\end{equation}
is a confidence interval for $Z$. Note that this second interval does
not depend on the initial composition of the urn, which could be
unknown.
\\[5pt]

\indent A remark useful for applications follows.

\begin{remark}
  \rm The estimators of $m,\, q_A,\, q_B$ and $q_{AB}$ defined in
  \eqref{eq-estimators} presuppose that we can observe both $A_j$ and
  $B_j$ for each $j=1,\dots, n$. Actually, in applications, we can
  observe $A_j$ (respectively, $B_j$) only when $X_j>0$ (respectively,
  $X_j<N_j$). Therefore, it makes more sense to use the following estimators:
  \begin{equation}\label{eq-estimators-v2}
    \begin{split}
  \widehat{m}_{n}&=\frac{\sum_{j=1}^n
    \left(A_jI_{\{X_j>0\}}+B_jI_{\{X_j=0\}}\right)}{n},\\
  \widehat{q}_{A,n}&=\frac{\sum_{j=1}^n
    A_j^2I_{\{X_j>0\}}}{\sum_{j=1}^nI_{\{X_j>0\}}}, \qquad
  \widehat{q}_{B,n}=\frac{\sum_{j=1}^n
    B_j^2I_{\{X_j<N_j\}}}{\sum_{j=1}^nI_{\{X_j<N_j\}}},\\ \qquad
  \widehat{q}_{AB,n}&=\frac{\sum_{j=1}^n
    A_jB_jI_{\{0<X_j<N_j\}}}{\sum_{j=1}^nI_{\{0<X_j<N_j\}}}\,.
  \end{split}
\end{equation}
    Note that $\widehat{m}_n\stackrel{a.s.}\longrightarrow m$ by Lemma
    \ref{lemma-app} (applied with
    $Y_j=A_jI_{\{X_j>0\}}+B_jI_{\{X_j=0\}}\leq C$ and
    $\mathcal{F}_j=\mathcal{G}_j$).  Indeed, we have
    \begin{equation*}
      \begin{split}
    E[A_jI_{\{X_j>0\}}+B_jI_{\{X_j=0\}}|\mathcal{G}_{j-1}]&= E\left[
      E[A_jI_{\{X_j>0\}}+B_jI_{\{X_j=0\}}|\mathcal{H}_{j-1}]\,
      |\mathcal{G}_{j-1}\right]\\ &=E[A_j
      P(X_j>0|\mathcal{H}_{j-1})+B_jP(X_j=0|\mathcal{H}_{j-1})
      \,|\mathcal{G}_{j-1}]=m_j\,,
      \end{split}
    \end{equation*}
    where the last equality is due to the fact that the conditional
    distribution of $X_j$ given $\mathcal{H}_{j-1}$ depends on $N_j,\,
    S_{j-1}$ and $H_{j-1}$ (and so coincides with the one given
    $\mathcal{G}_{j-1}$) and to relation \eqref{sp-cond-rinforzi}. The
    convergence $\widehat{q}_{A,n}\stackrel{a.s.}\longrightarrow q_A$
    also follows from by Lemma \ref{lemma-app}. Indeed, we have
    $$
    E[A^2_jI_{\{X_j>0\}}|\mathcal{H}_{j-1}]
    =A_j^2\left[1-\frac{\binom{S_{j-1}-H_{j-1}}{N_j}}{\binom{S_{j-1}}{N_j}}\right]
    $$ Then, conditioning with respect to $\mathcal{G}_{j-1}$ and
    using \eqref{sp-cond-rinforzi}, we get
    $E[A^2_jI_{\{X_j>0\}}|\mathcal{G}_{j-1}]=q_{A,j}\varphi(N_j,S_{j-1},H_{j-1})$
    with $\varphi(N,S,H)=
    \left[1-\frac{\binom{S-H}{N}}{\binom{S}{N}}\right]$. Finally,
    conditioning with respect to $\mathcal{F}_{j-1}$, we find
    $$
    E[A^2_jI_{\{X_j>0\}}|\mathcal{F}_{j-1}]=q_{A,j}\sum_{k=1}^C
    \varphi(k,S_{j-1},H_{j-1})P(N_j=k|\mathcal{F}_{j-1})\,.
    $$ Assuming that
    $P(N_j=k|\mathcal{F}_{j-1})\stackrel{a.s.}\longrightarrow \nu(k)$
    (with $\nu(k)$ possibly random), as a consequence of Proposition
    \ref{conv-S} and the above equality, we have
    $$
    E[A^2_jI_{\{X_j>0\}}|\mathcal{F}_{j-1}]\stackrel{a.s.}\sim
    q_{A,j}\sum_{k=1}^C
    \left[1-(1-\frac{H_{j-1}}{S_{j-1}})^k\right]P(N_j=k|\mathcal{F}_{j-1})
    \stackrel{a.s.}\longrightarrow
    q_A\sum_{k=1}^C
    \left[1-(1-Z)^k\right]\nu(k)\,.
    $$ Similarly, we have
    $E[I_{\{X_j>0\}}|\mathcal{F}_{j-1}]\stackrel{a.s.}\longrightarrow
    \sum_{k=1}^C \left[1-(1-Z)^k\right]\nu(k)$ and so, by Lemma
    \ref{lemma-app}, we obtain
    $$
\widehat{q}_{A,n}=\frac{\sum_{j=1}^n
  A_j^2I_{\{X_j>0\}}/n}{\sum_{j=1}^nI_{\{X_j>0\}}/n}\stackrel{a.s.}
\longrightarrow
\frac{q_A \sum_{k=1}^C
    [1-(1-Z)^k]\nu(k)}{\sum_{k=1}^C
    [1-(1-Z)^k]\nu(k)}=q_A\,.
$$ Exactly with the same argument, we get
$\widehat{q}_{B,n}\stackrel{a.s.}\longrightarrow q_B$.  For the almost
sure convergence of $\widehat{q}_{AB,n}$ to $q_{AB}$, we can argue in
the similar way, but we need $P(\nu(1)<1)=1$ in order to guarantee
that $\sum_{k=1}^C \left[1-(1-Z)^k-Z^k\right]\nu(k)>0$ almost surely.
\end{remark}

\section{Examples and numerical illustrations}
\label{sec:simulations}

Before considering special cases as illustration through numerical
simulations, let us formulate some general remarks.

\begin{remark}\label{rem:A-B}{\em ($[A_n,B_n]$ identically distributed)} \rm 
If all the random vectors $[A_n,B_n]$ (that are independent by
assumption~(A3)) are also identically distributed, then we simply have
$m=m_n=E[A_n]=E[B_n]$ and condition~\eqref{ass-A-B} is satisfied with
$q_A=q_{A,n}=E[A_n^2],\,q_B=q_{B,n}=E[B_n^2]$ and
$q_{AB}=q_{AB,n}=E[A_nB_n]$).
\end{remark}

\begin{remark} \label{rem:4.2} {\em ($N_n$ independent of the past)}\\
If, for each $n$, the random variable $N_n$ is independent of
${\mathcal F}_{n-1}$, then we simply have
$E[N_n|\mathcal{F}_{n-1}]=E[N_n]$,
$E[N_{n}^2|\mathcal{F}_{n-1}]=E[N_n^2]$ and
$E[N_n^{-1}|\mathcal{F}_{n-1}]=E[N_n^{-1}]$.  Therefore, conditions
\eqref{ass-base-2}, \eqref{ass-N-sp-cond-2} and~\eqref{ass-base-bis}
are satisfied whenever the above sequences of mean values converge to
suitably constants~$N$, $Q$ and $L$. For instance, this happens when
all the random variables~$N_n$ are identically di\-stri\-bu\-ted.
More precisely, in this last case, assuming $N_n\leq a+b$ (so that we
are sure that $N_{n}\leq S_{n-1}$ for each $n$), with mean value $\mu$
and variance $\sigma^2$, conditions~\eqref{ass-base-2},
\eqref{ass-N-sp-cond-2} and~\eqref{ass-base-bis} are satisfied with
$N=E[N_n]=\mu$, $Q=E[N_n^2]=q_N=\sigma^2+\mu^2$ and
$L=E[N_n^{-1}]=\eta$.
\end{remark}

\begin{remark}\label{rem-ipotesi} {\em ($N_n$ dependent on $Z_{n-1}$)}
\rm When $N_n$ depends on the urn proportion at time $n-1$,
i.e.~$Z_{n-1}$, in such a way that, for each $n$, we have
  $$ E[N_{n+1}|\mathcal{F}_{n}]=f(Z_{n}), \quad 
  E[N_{n+1}^2|\mathcal{F}_{n}]=g(Z_{n}), \quad  
  E[N_{n+1}^{-1}|\mathcal{F}_{n}]=h(Z_{n})\,,
  $$ where $f,\,g$ and $h$ are continuous functions, then
  conditions~\eqref{ass-base-2}, \eqref{ass-N-sp-cond-2}
  and~\eqref{ass-base-bis} are satisfied with $N=f(Z)$, $Q=g(Z)$ and
  $L=h(Z)$. Note that, if the functions $f,\,g$ and $h$ are known, we
  can obtain asymptotic confidence intervals for $Z$ replacing
  $\widehat{\mu}_n$ and $\widehat{q}_{N,n}$ in the expression for
  $V_n$ by $f(Z_n)$ and $g(Z_n)$, respectively, and replacing
  $\widehat{\mu}_n,\,\widehat{q}_{N,n}$ and $\widehat{\eta}_n$ in the
  expression for $W_n$ by $f(M_n),\,g(M_n)$ and $h(M_n)$,
  respectively.
\end{remark}   

\begin{remark}\label{remark-example}
  \label{rem:4.3}{\em ($N_n$ almost surely convergent)}\\
  \rm If $(N_n)$ is a sequence of integer-valued random variables with
  $1\leq N_n\leq C$ and converging almost surely to a random variable
  $N$, then (by Lemma~\ref{lemma-app-C}) conditions
  \eqref{ass-base-2}, \eqref{ass-N-sp-cond-2} and \eqref{ass-base-bis}
  are satisfied and $Q=N^2$ and $L=N^{-1}$. See, for instance,
  Example 4.2 in~\cite{crimaldi2016}, where $(N_n)$ is a symmetric
  random walk with two absorbing barriers.
\end{remark}
\medskip

\indent The following examples regard the case 2) (that is the case of
equal reinforcement means) and they deal with the different situations
described in the above general remarks.

\begin{figure}
\centering
\begin{subfigure}[t]{0.49\linewidth}
\centering
\includegraphics[scale=0.85]{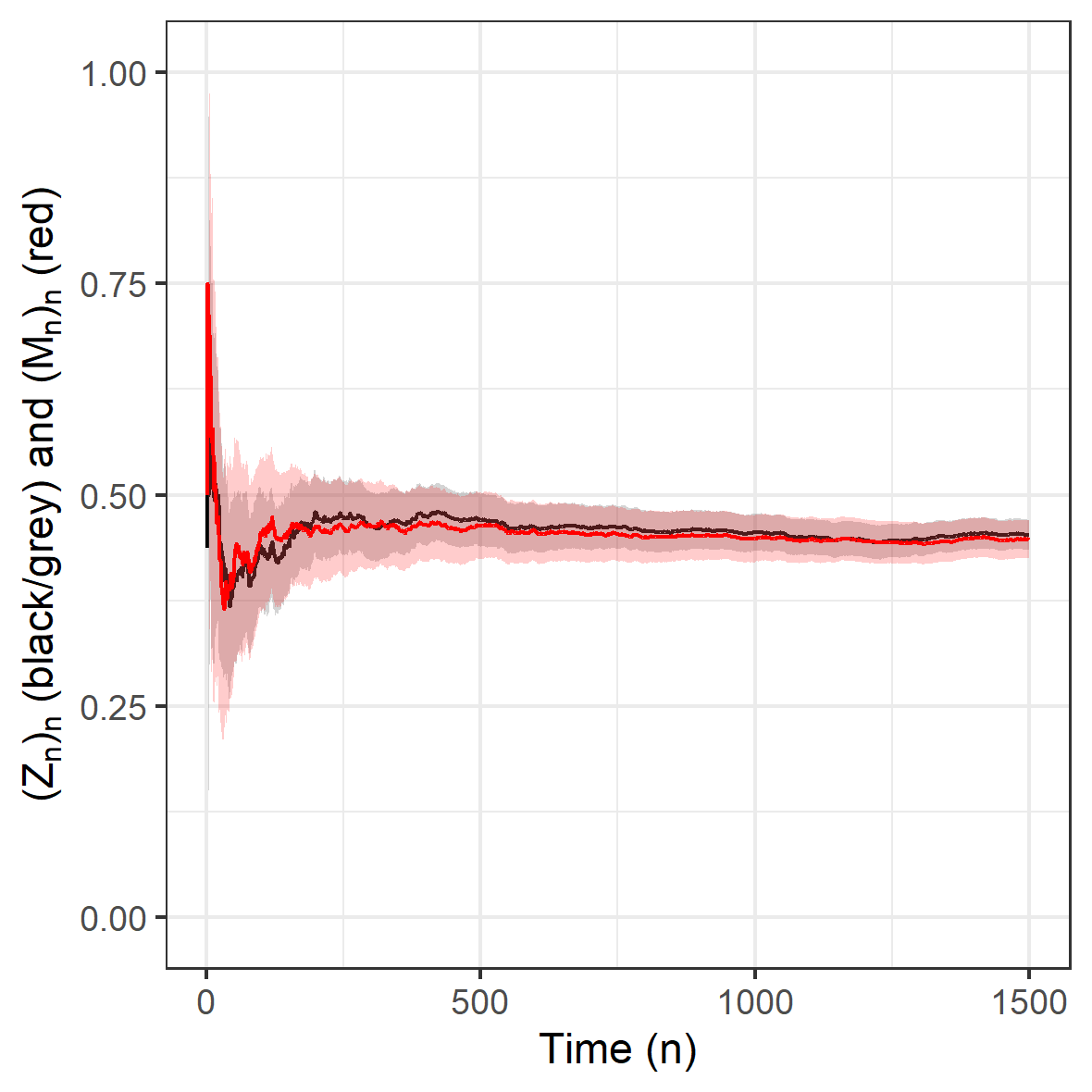}
\end{subfigure}
\hfill 
\begin{subfigure}[t]{0.49\linewidth}
\centering
\includegraphics[scale=0.85]{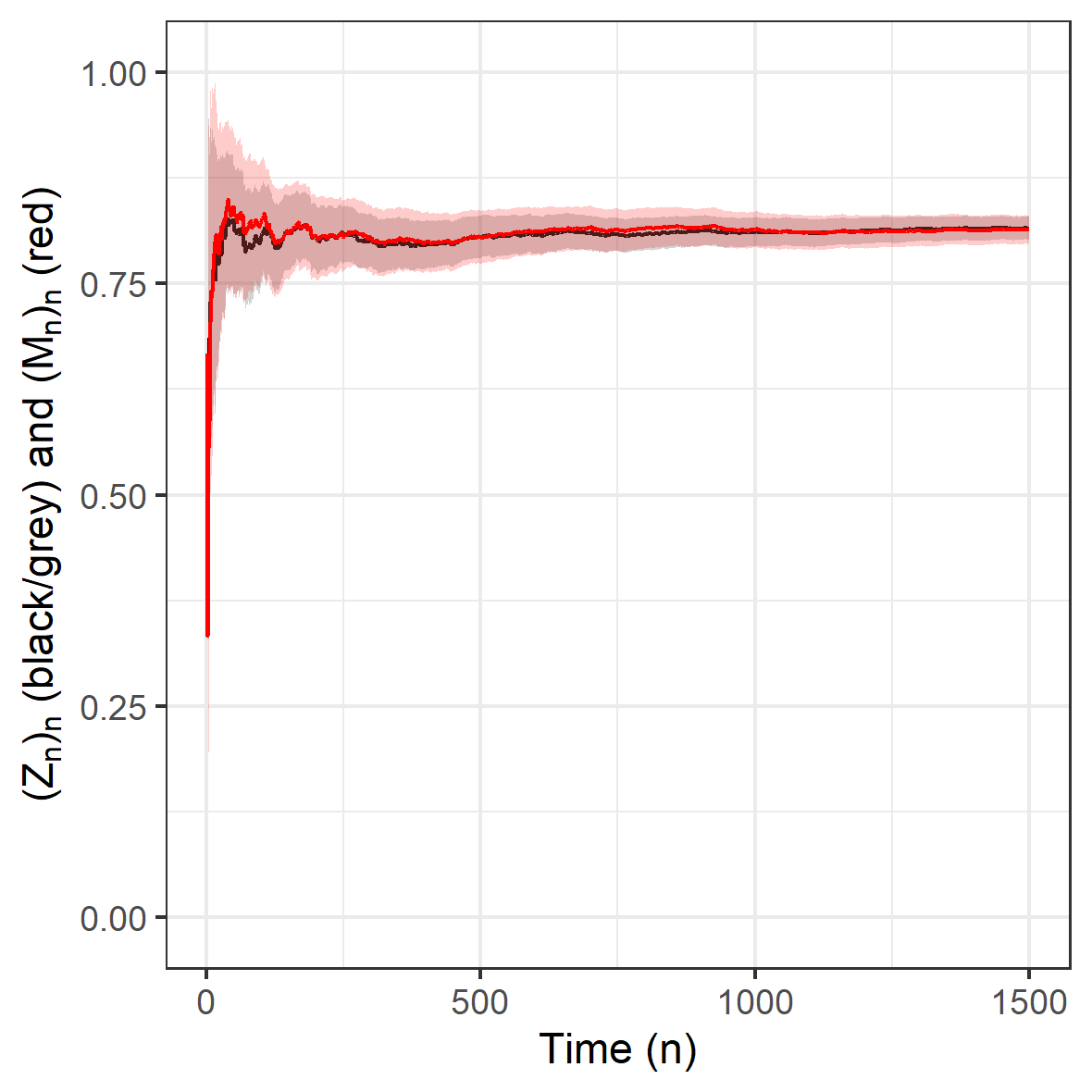}
\end{subfigure}
\caption{Case 1a. Time-horizon $1500$.  On each picture, one sample
  plot of $(Z_n)_n$ (black) and $(M_n)_n$ (red) with the corresponding
  confidence intervals for~$Z$ with $\alpha=0.05$ (resp. grey and
  red).}
\label{fig:example1a}
\end{figure}

\noindent \textbf{Example 1a} \\ Take each $N_n$ independent of
$\mathcal{F}_{n-1}$ and uniformly distributed on
$\{1,\hdots,5\}$. Moreover, take $A_n$ and $B_n$ satisfying assumption
(A3), independent and uniformly distributed on $\{1,\hdots,5\}$. We
set $a=b=5$. See Fig.~\ref{fig:example1a} for samples.  \\

\begin{figure}
\centering
\begin{subfigure}[t]{0.49\linewidth}
\centering
\includegraphics[scale=0.85]{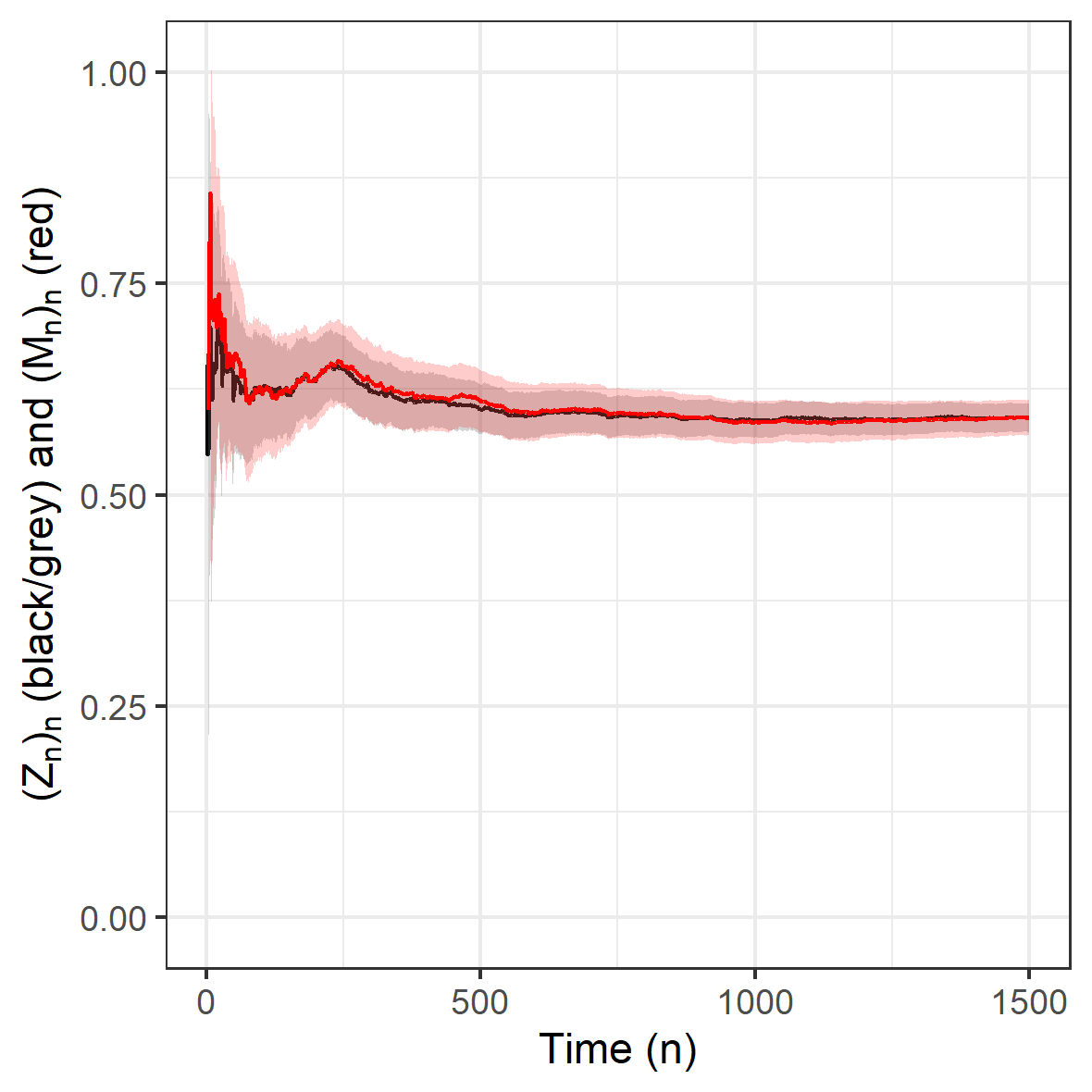}
\end{subfigure}
\hfill 
\begin{subfigure}[t]{0.49\linewidth}
\centering
\includegraphics[scale=0.85]{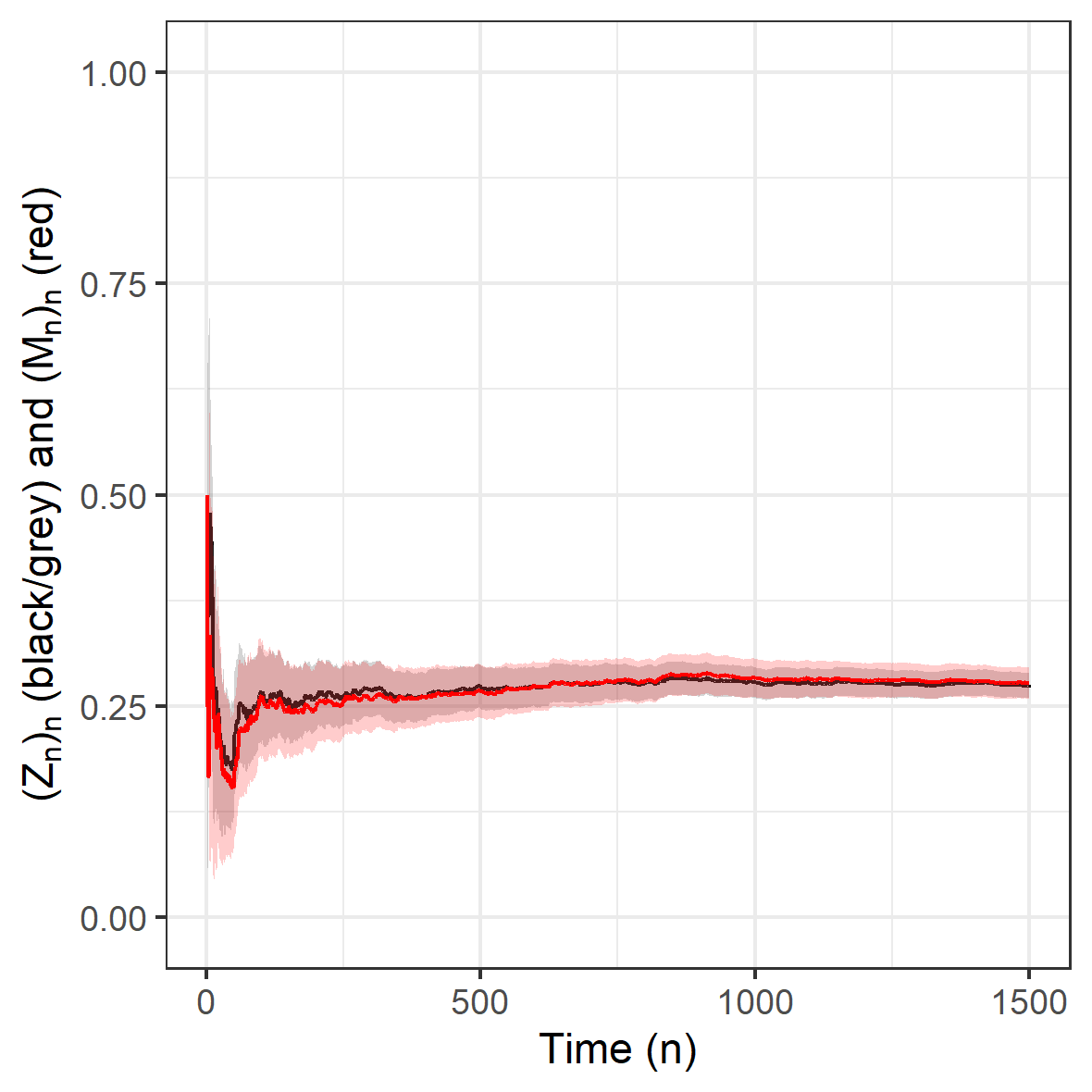}
\end{subfigure}
\caption{Case 1b. Time-horizon $1500$. On each picture, one sample
  plot of $(Z_n)_n$ (black) and $(M_n)_n$ (red) with the corresponding
  confidence intervals for~$Z$ with $\alpha=0.05$ (resp. grey and
  red).}
\label{fig:example1b}
 \end{figure}

 \noindent \textbf{Example 1b} \\ Take each $N_n$ independent of
 $\mathcal{F}_{n-1}$ and uniformly distributed on $\{1,\hdots,5\}$.
 In particular, assumption (i) in Section~\ref{conf-int} is satisfied.
 Moreover, take $[A_n,B_n]$ satisfying assumption~(A3) and such that
 $$
A_n\stackrel{d}=1+Y_1\qquad\mbox{and}\qquad B_n\stackrel{d}=1+Y_2\,,
$$
where $Y_1$ and $Y_2$ are, respectively,  the first and the second
 component of a multinomial distribution associated to the parameters:
 size$=12$, probabilities$=(4/15,4/15,7/15)$. Thus the random
 variables $A_n$ and $B_n$ are negatively correlated. We set $a=b=5$.
 See Fig.~\ref{fig:example1b} for samples.  \\

\begin{figure}
\centering
\begin{subfigure}[t]{0.49\linewidth}
\centering
\includegraphics[scale=0.85]{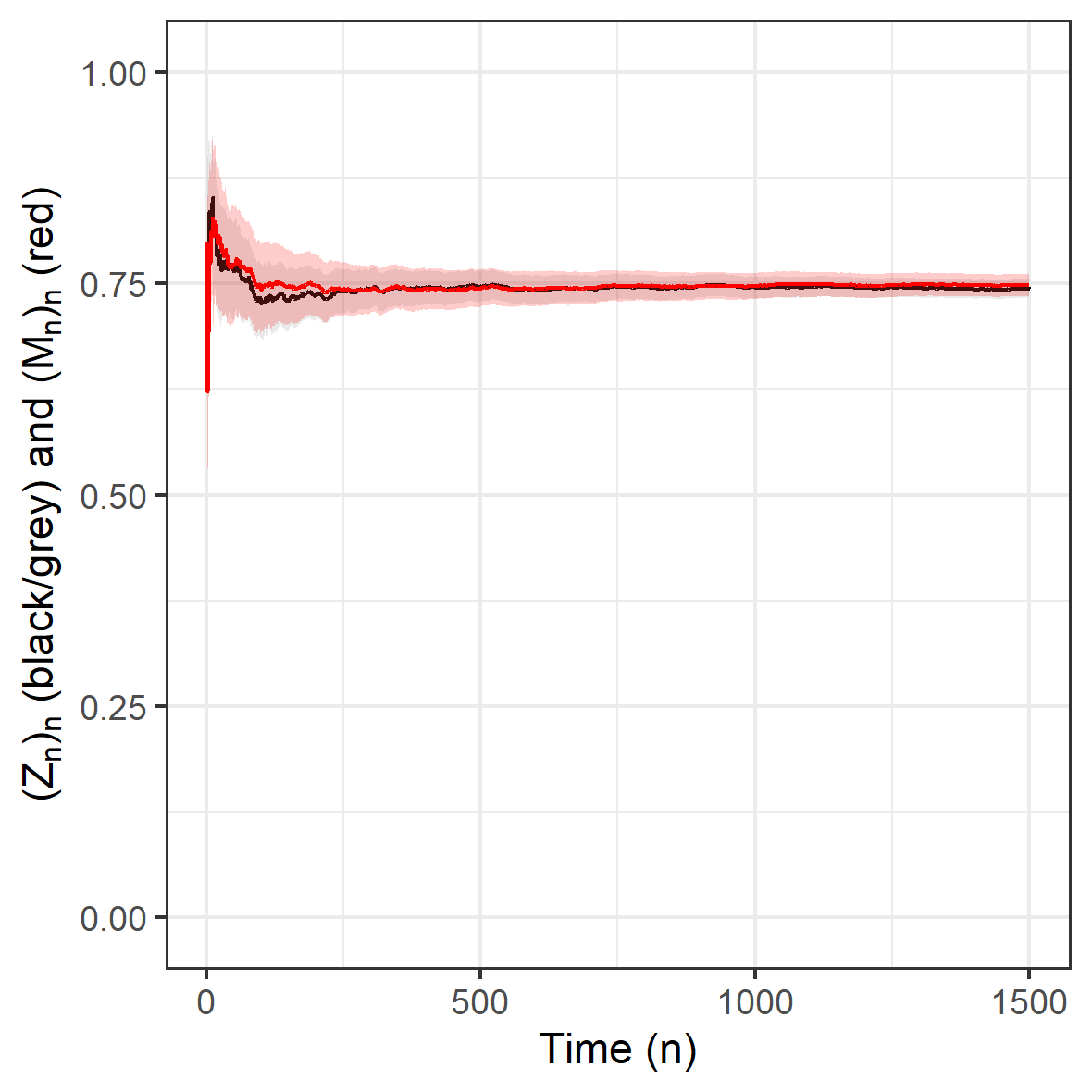}
\end{subfigure}
\hfill 
\begin{subfigure}[t]{0.49\linewidth}
\centering
\includegraphics[scale=0.85]{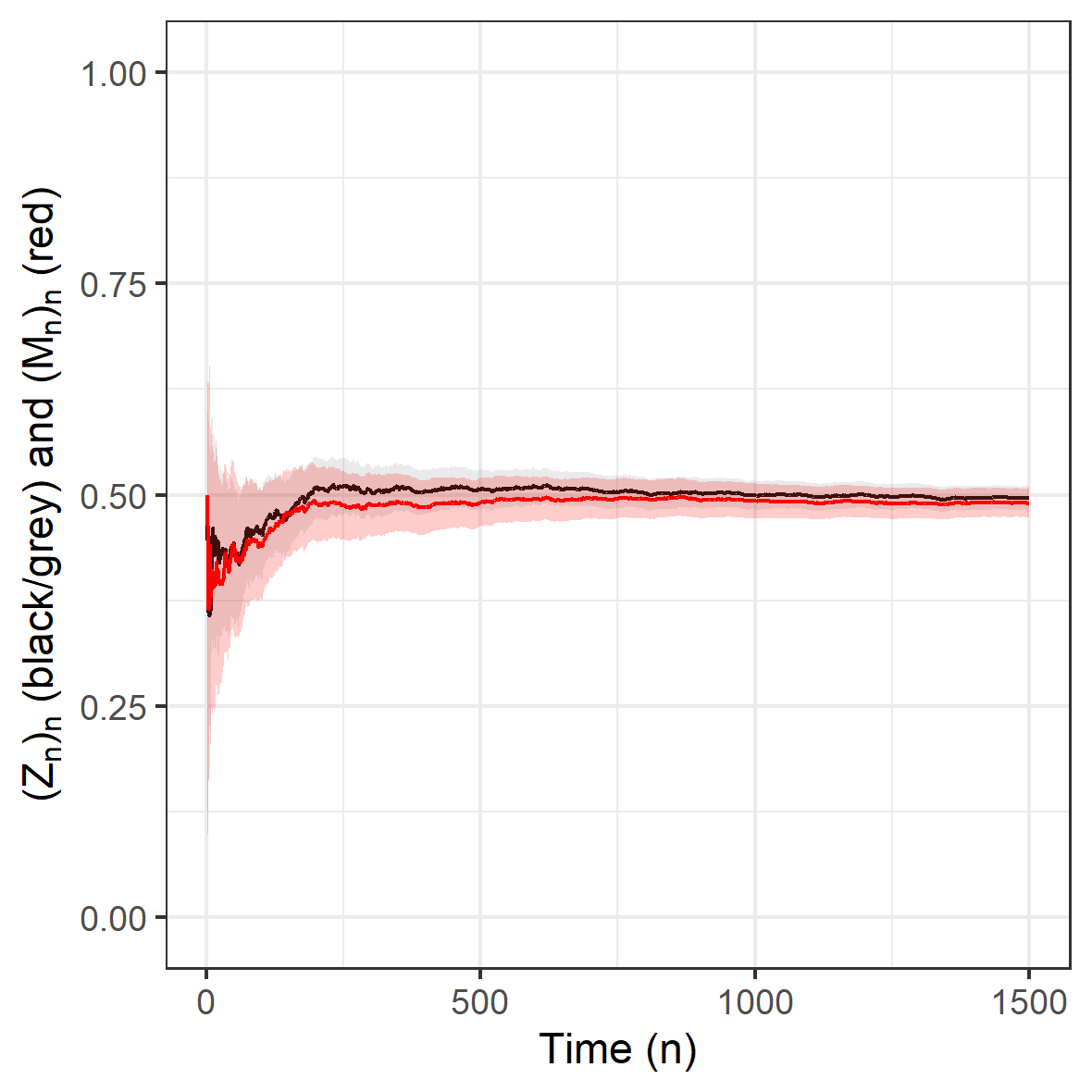}
\end{subfigure}
\caption{Case 1c. Time-horizon $1500$. On each picture, one sample
  plot of $(Z_n)_n$ (black) and $(M_n)_n$ (red) with the corresponding
  confidence intervals for~$Z$ with $\alpha=0.05$ (resp. grey and
  red). The confidence intervals are the ones given in Remark
  \ref{rem-ipotesi}, taking parameter $\kappa$ known (that is with the
  functions $f,\,g$ and $h$ known).}
\label{fig:example1c}
\end{figure}

 \noindent \textbf{Example 1c} \\ Set $(N_n)_n$ be a
 sequence of random variables such that
$$N_n | {\mathcal F}_{n-1} \stackrel{d}= 1+ {\mathcal
   B}(\kappa,Z_{n-1}).$$ Moreover, take $A_n$ and $B_n$ satisfying
 assumption (A3), independent and uniformly distributed on
 $\{1,\hdots,5\}$. In particular, we are in the situation described in
 Remark \ref{rem-ipotesi}. Indeed, we have:
\begin{equation*}
  \begin{split}
    E[N_{n+1}|\mathcal{F}_n]&=1+\kappa Z_n\stackrel{a.s.}\longrightarrow
    N=1+ \kappa Z
\\
E[N_{n+1}^2|\mathcal{F}_n]&=\kappa Z_n(1-Z_n)+(1+\kappa Z_n)^2
\stackrel{a.s.}\longrightarrow Q=\kappa Z(1-Z)+(1+ \kappa Z)^2
\\
E[N_{n+1}^{-1}|\mathcal{F}_n]&=\frac{1-(1-Z_n)^{\kappa +1}}{(\kappa +1) Z_n}
\stackrel{a.s.}\longrightarrow \frac{1-(1-Z)^{\kappa +1}}{(\kappa +1) Z}
\end{split}
  \end{equation*}
  (recall that $P(Z=0)=0$ by Theorem \ref{th-legge-2}).  We set
$\kappa=10$ and $a=b=6$. See Fig.~\ref{fig:example1c} for samples.\\

%
%
%
%
%

\begin{figure}
\centering
\begin{subfigure}[t]{0.49\linewidth}
\centering
\includegraphics[scale=0.85]{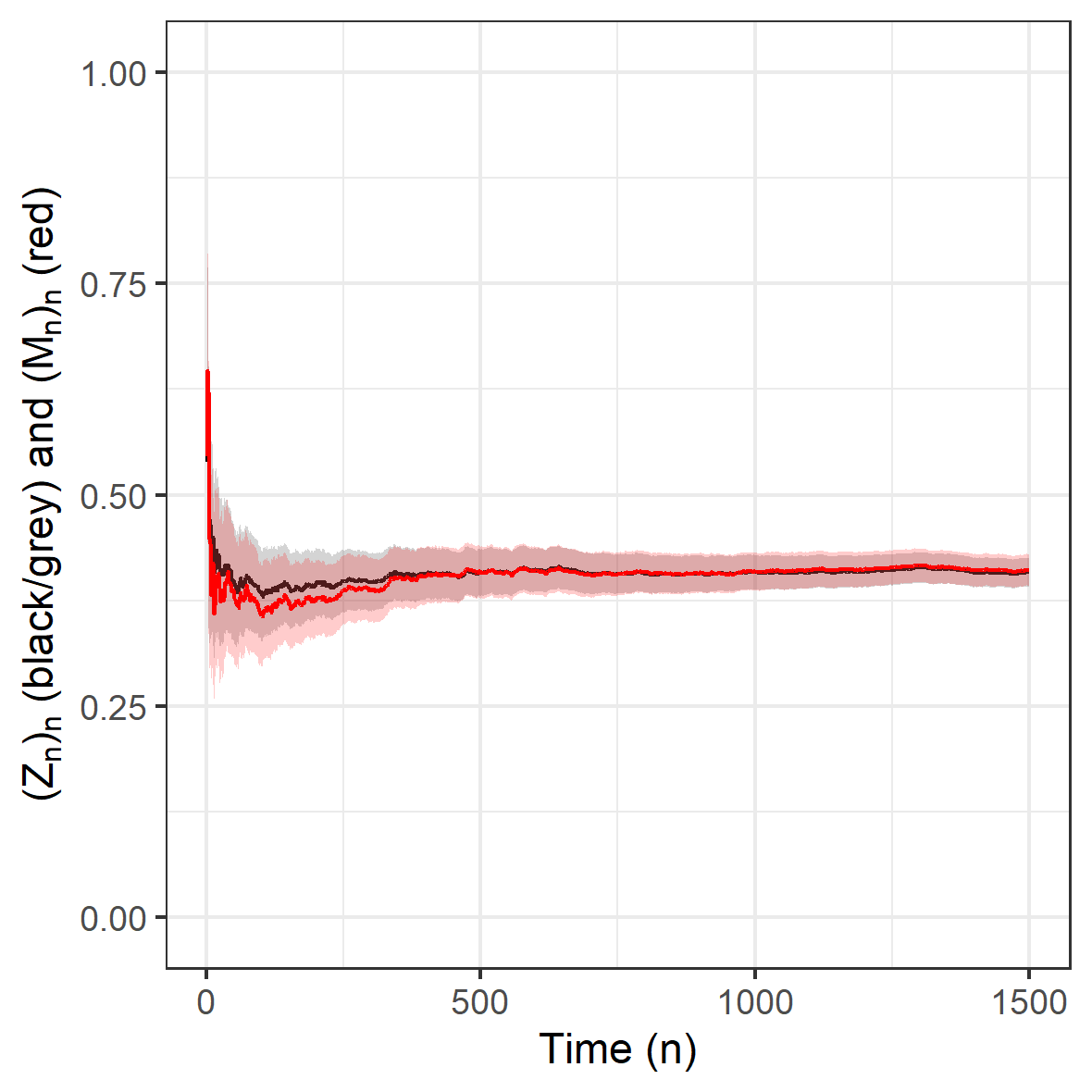}
\end{subfigure}
\hfill 
\begin{subfigure}[t]{0.49\linewidth}
\centering
\includegraphics[scale=0.85]{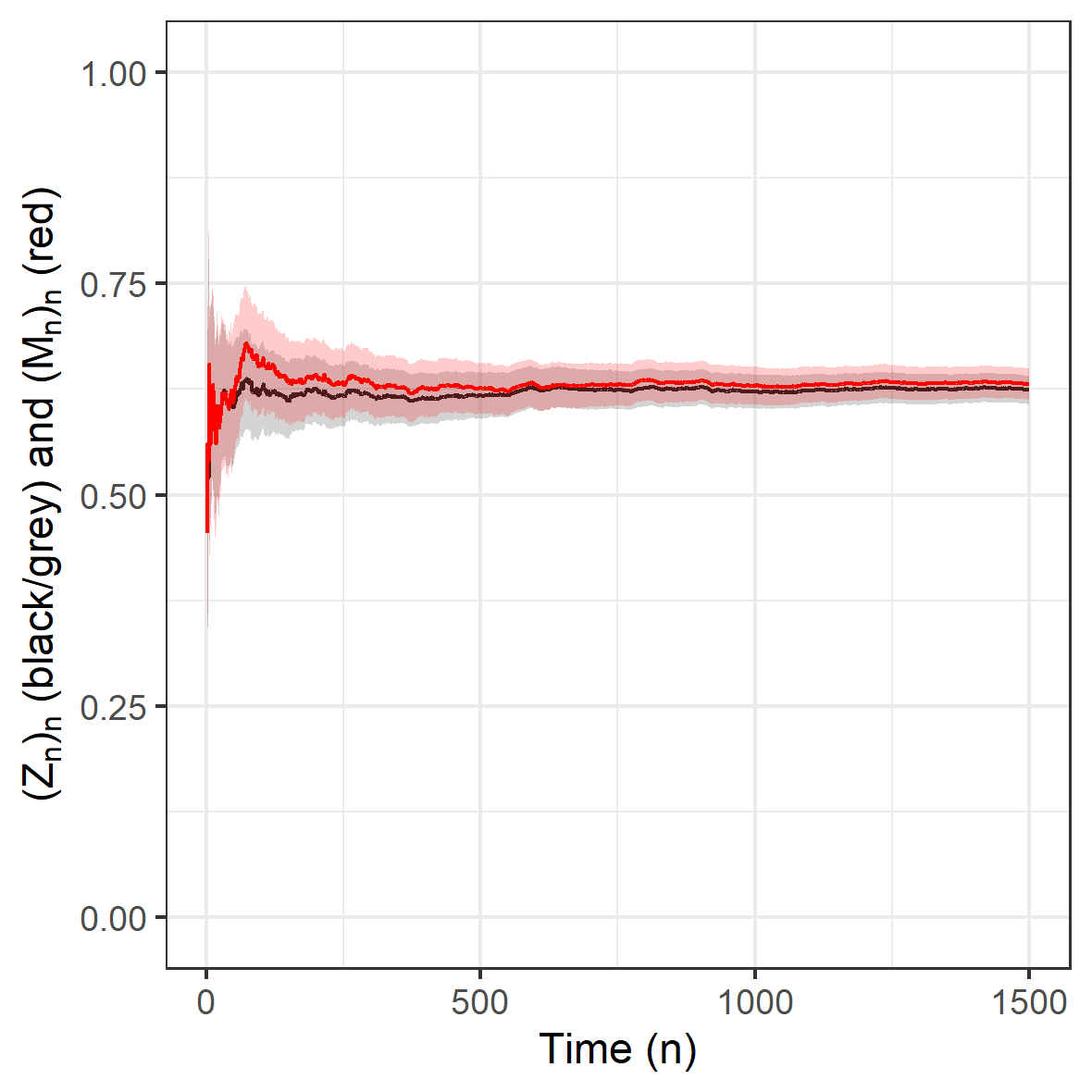}
\end{subfigure}
\caption{Case 1d.  Time-horizon $1500$.   On each picture, one sample
 plot of $(Z_n)_n$ and $(M_n)_n$  with the corresponding confidence
intervals for~$Z$ with $\alpha=0.05$ (resp. grey and red). }
\label{fig:example2a}
\end{figure}

 \noindent \textbf{Example 1d}\\
 %
 %
 %
Take each $N_n$ independent of $\mathcal{F}_{n-1}$ and such that
$$ N_n \stackrel{d}= 2+\mathcal B(\kappa,p_n),$$ with $\kappa=10$ and
$p_n=1/\sqrt{n}$. Moreover, take $A_n$ and $B_n$ satisfying assumption
(A3), independent and such that
$$ A_n \stackrel{d}= B_n \stackrel{d}= 1+\mathcal B(\kappa',q_n),$$
with $\kappa'=5$ and $q_n=\min(1,\frac12 + \frac{1}{\sqrt{n}})$. We
take $a=b=6$.  See Fig.~\ref{fig:example2a} for samples.  \\

\begin{figure}
\centering
\begin{subfigure}[t]{0.49\linewidth}
  \centering
\includegraphics[scale=0.85]{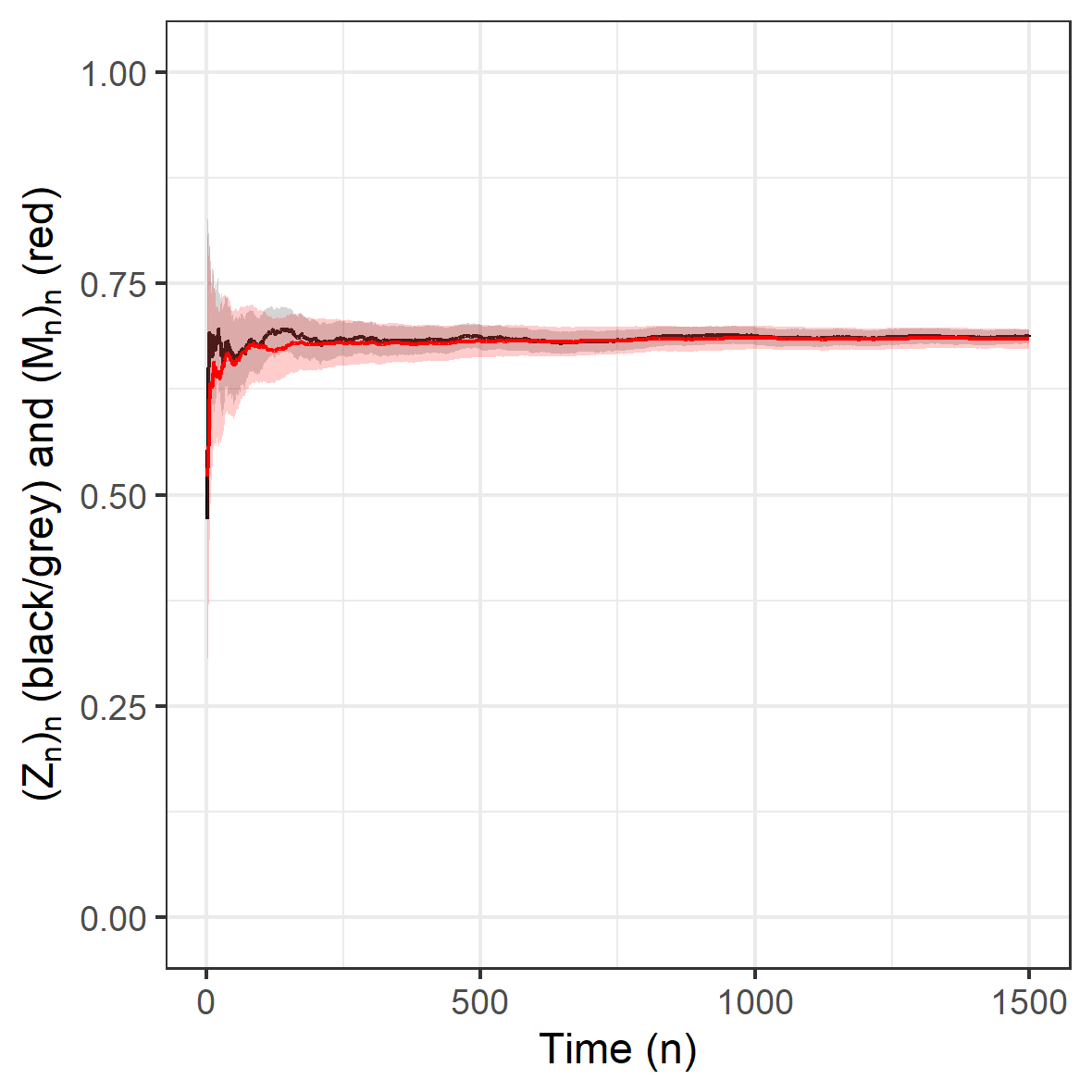}
\end{subfigure}
\hfill 
\begin{subfigure}[t]{0.49\linewidth}
  \centering
  \includegraphics[scale=0.85]{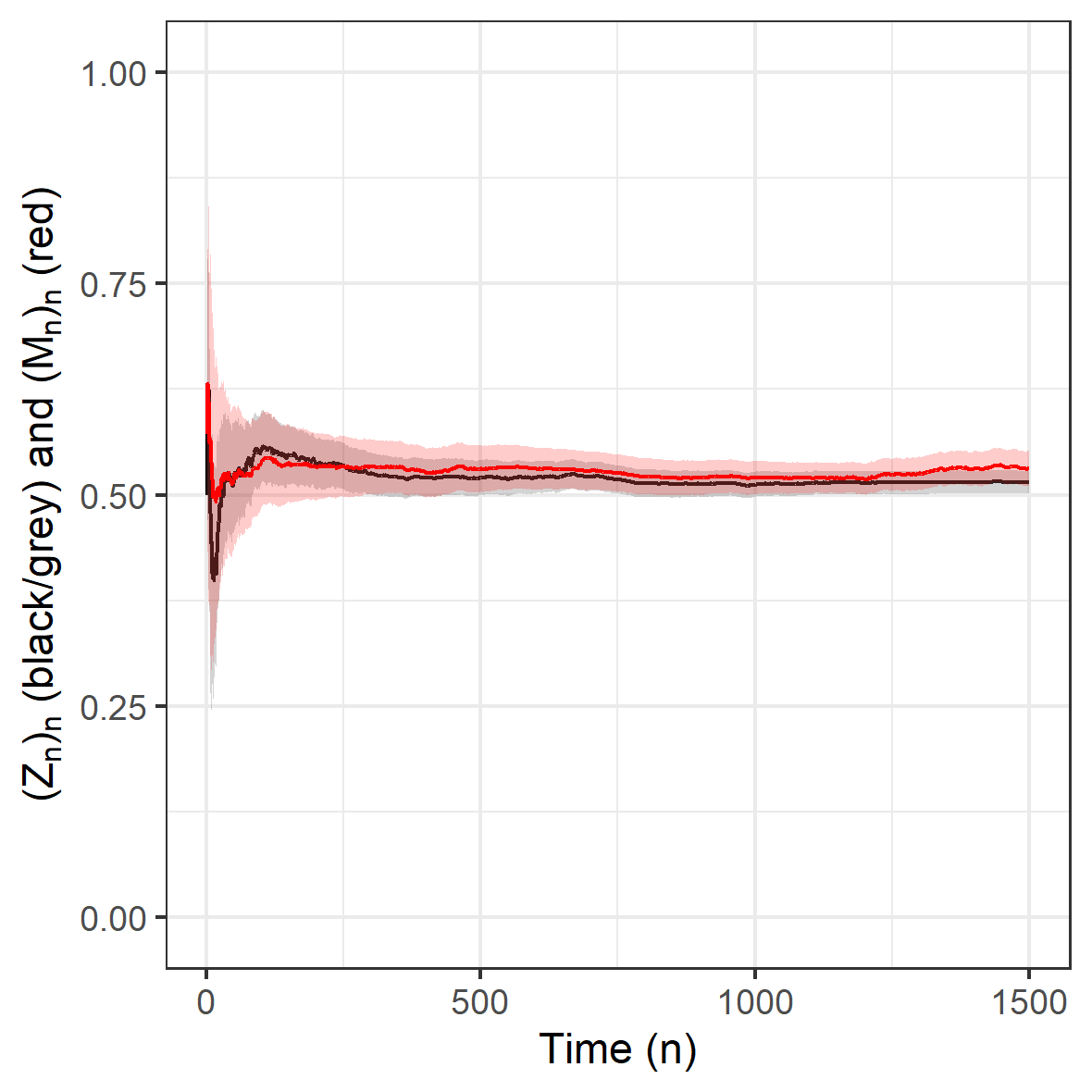}
\end{subfigure}
\caption{Case 1e. Time-horizon $1500$.  On each picture, one sample
  plot of $(Z_n)_n$ and $(M_n)_n$ with the corresponding confidence
  intervals for~$Z$ with $\alpha=0.05$ (resp. grey and red).}
\label{fig:example2b}
\end{figure}

 \noindent \textbf{Example 1e} \\
%
%
This example is associated to Remark~\ref{rem:4.3}. Following
 Example~4.2 in~\cite{crimaldi2016}, take $(N_n)_n$ be a sequence of
 random variables defined through a symmetric nearest neighbors random
 walk with absorbing barriers.  Given $h\in{\mathbb N}$, with $3\leq h
 \leq a+b$, let $\widetilde{N}_1$ be a random variable with
 distribution concentrated on $\{2,\dots, h-1\}$ and set
$$
\widetilde{N}_n=\widetilde{N}_1+\sum_{j=2}^{n} Y_j \textrm{ for }
n \geq 2\,,$$
$$
T_1=\inf\{n:\widetilde{N}_n=1\},\qquad
T_h =\inf\{n:\widetilde{N}_n=h\}
$$
and 
$$ N_n=\widetilde{N}_{T\wedge n}\quad\mbox{for } n\geq 1, \quad 
\mbox{with } T=T_1\wedge T_h\,,
$$ where each $Y_j$ is independent of
$[\widetilde{N}_1,X_1,A_1,B_1,Y_1, X_2, A_2,B_2,\dots, Y_{j-1},
  X_{j-1}, A_{j-1},B_{j-1}]$ and such that
$P(Y_j=-1)=P(Y_j=1)=p \in
(0,\frac12]$ and $P(Y_j=0)=1-2p$.  Then
  $N_n\stackrel{a.s.}\longrightarrow N=\widetilde{N}_T$ where
  $N=\mathds{1}_{\{T=T_1\}}+h\mathds{1}_{\{T=T_h\}}$. We take
  $A_n$ and $B_n$ satisfying assumption~(A3), independent and such that
$$ A_n \stackrel{d}= B_n \stackrel{d}= 1+\mathcal B(\kappa',q_n).$$ We
  consider specifically $a=b=30$, $h=50$, $\widetilde{N}_1$ uniformly
  distributed on $\{2,\dots,h-1\}$, $p=1/4$, $\kappa'=5$ and
  $q_n=\min(1,\frac12 + \frac{1}{\sqrt{n}})$.  Note that also in this
  case it is possible to contruct confidence intervals for $Z$ (see
  Remark \ref{rem-ipotesi}). See Fig.~\ref{fig:example2b} for samples.
  \\

\noindent In the following example, the random variables $N_n$, $A_n$
and $B_n$ are not bounded, but condition~\eqref{cond-quasi-mart} is
satisfied.
\\

\begin{figure}
\centering
\begin{subfigure}[t]{0.49\linewidth}
\centering
\includegraphics[scale=0.85]{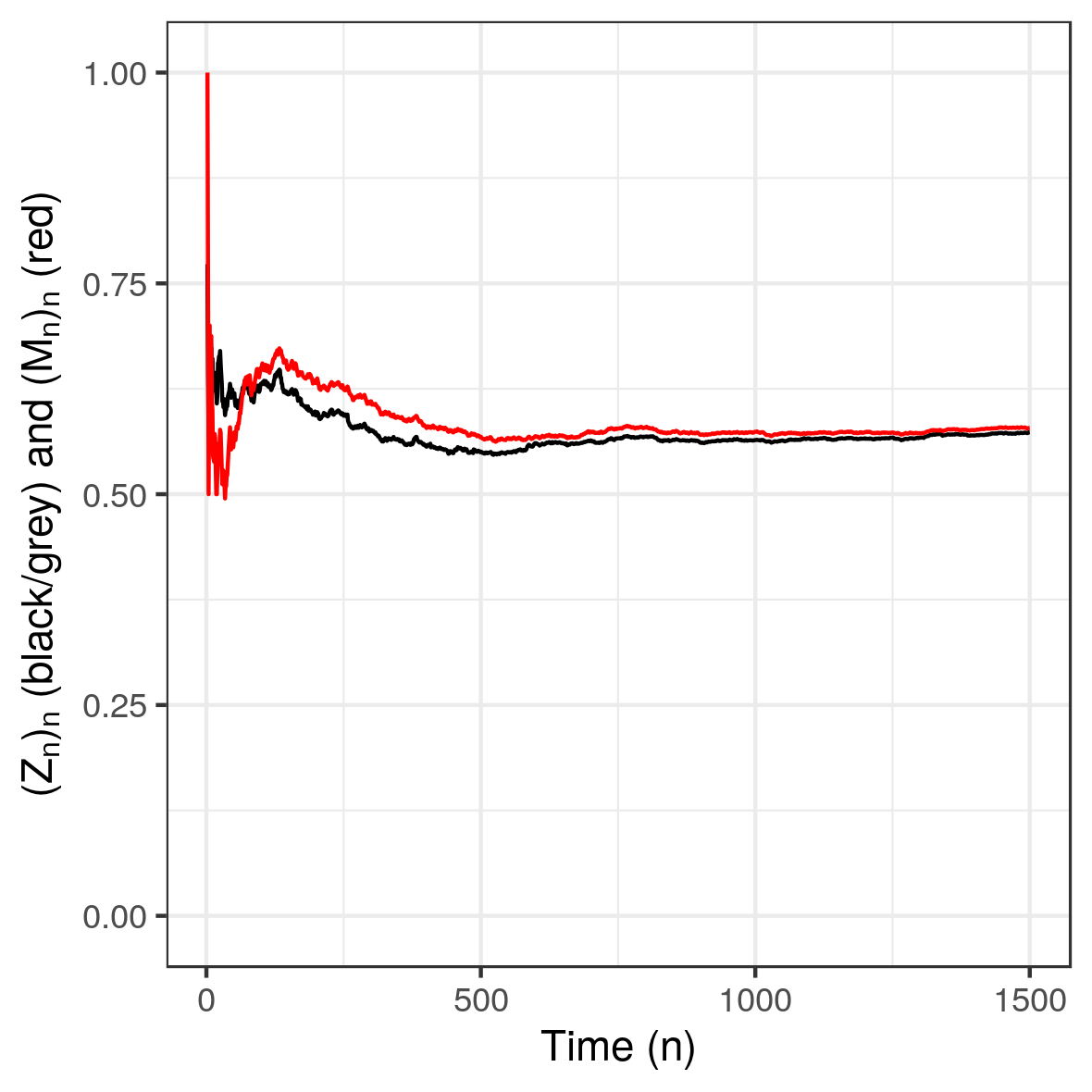}
\end{subfigure}
\hfill 
\begin{subfigure}[t]{0.49\linewidth}
\centering
\includegraphics[scale=0.85]{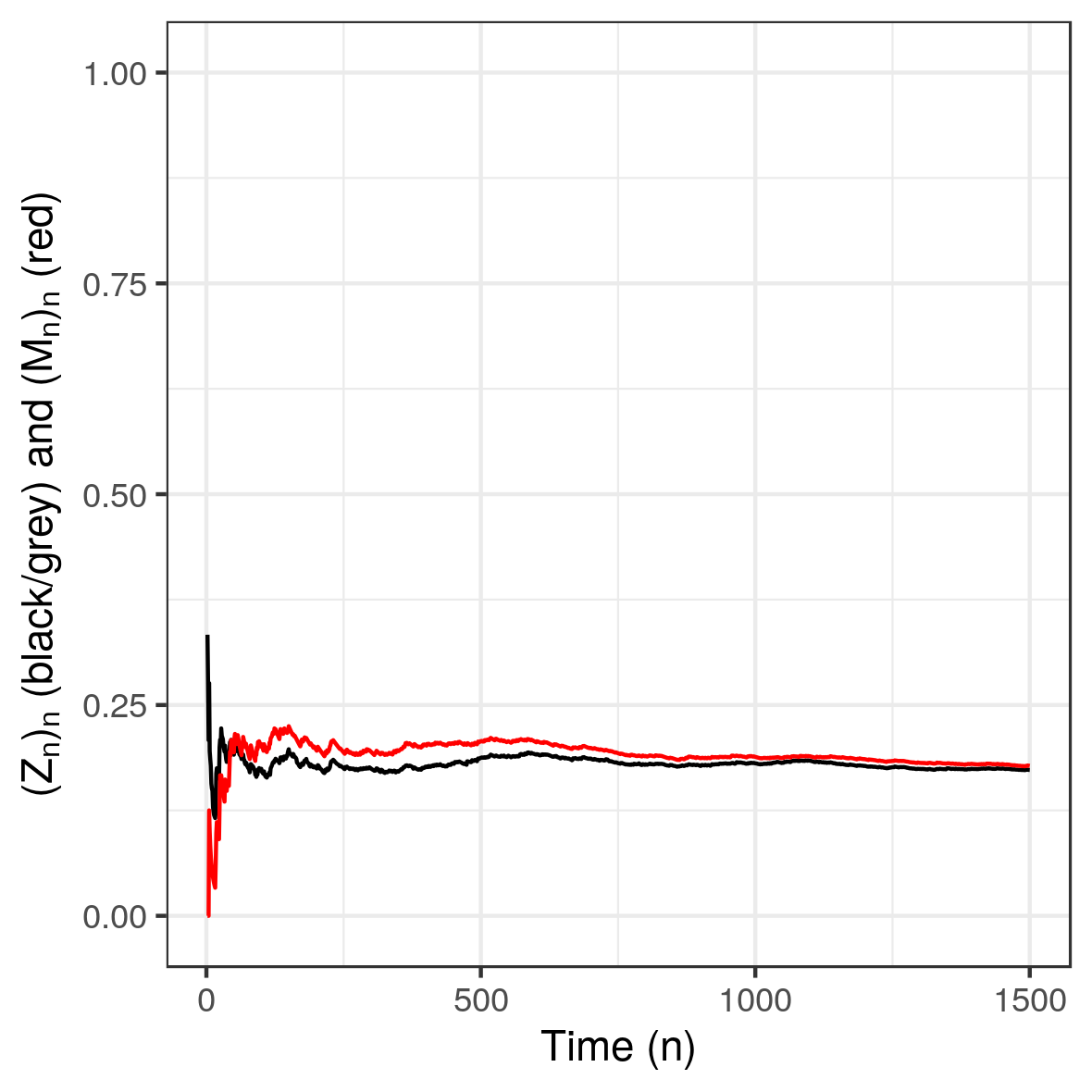}
\end{subfigure}
\caption{Case 2. Time-horizon $1500$.  On each picture, one sample
  plot of $(Z_n)_n$ and $(M_n)_n$ (resp. black and red).
  %
  %
}
\label{fig:example3a}
\end{figure}

\noindent \textbf{Example 2}\\
%
%
For each $n\geq 1$, take
 $\widetilde{N}_n$ independent of
 $[\widetilde{N}_1,X_1,A_1,B_1,\dots,\widetilde{N}_{n-1},X_{n-1},A_{n-1},B_{n-1}]$
 and such that
$$\widetilde{N}_n \stackrel{d}= 1 + {\mathcal
   B}(\kappa+\ceil*{n^{\frac13}},p)$$ with $\kappa=3$ and
 $p=1/10$. Set $N_n=\widetilde{N}_n\wedge S_{n-1}$ for each $n\geq 1$.
 Take $A_n$ and $B_n$ satisfying assumption~(A3), independent and such
 that
$$ A_n \stackrel{d}= B_n \stackrel{d}= 1 + \textrm{neg}{\mathcal
   B}(r,p_n)\,,$$ where $\textrm{neg}{\mathcal B}(r,p_n)$ means the
 negative binomial distribution with parameters $r=3$ and
 $p_n=1/\sqrt{n+1}$, that is with mean value equal to $rp_n/(1-p_n)$ and
 variance equal to $rp_n/(1-p_n)^2$. 
 Condition~\eqref{cond-quasi-mart} is satisfied because
$$
E[(A_n+B_n)^2]=O(1)\quad\mbox{and}\qquad
E[N_n^2]\leq E[\widetilde{N}_n^2]=O(n^{2/3}).
$$
We set $a=b=5$.  See
Fig.~\ref{fig:example3a} for samples.
\\

\indent The last two examples below are related to the case
$m_A>m_B$. Note that the time of the almost sure convergence to~$1$,
proven above, depends on the difference $m_A-m_B$. Thus, when this
difference is small, it may be difficult to guess the right asymptotic
behavior only through simulations.  \\

\begin{figure}
\centering
\begin{subfigure}[t]{0.49\linewidth}
\centering
\includegraphics[scale=0.85]{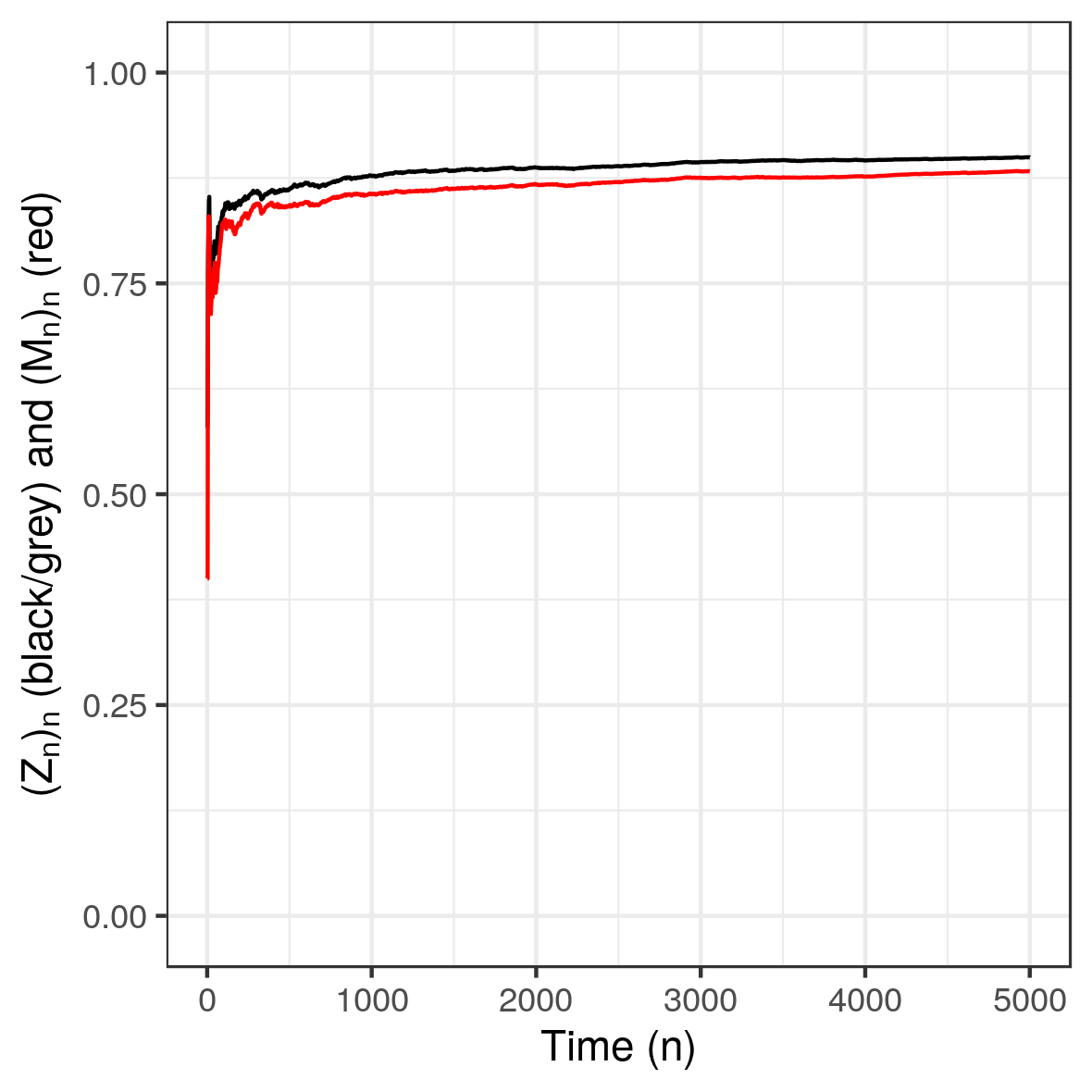}
\end{subfigure}
\hfill 
\begin{subfigure}[t]{0.49\linewidth}
\centering
\includegraphics[scale=0.85]{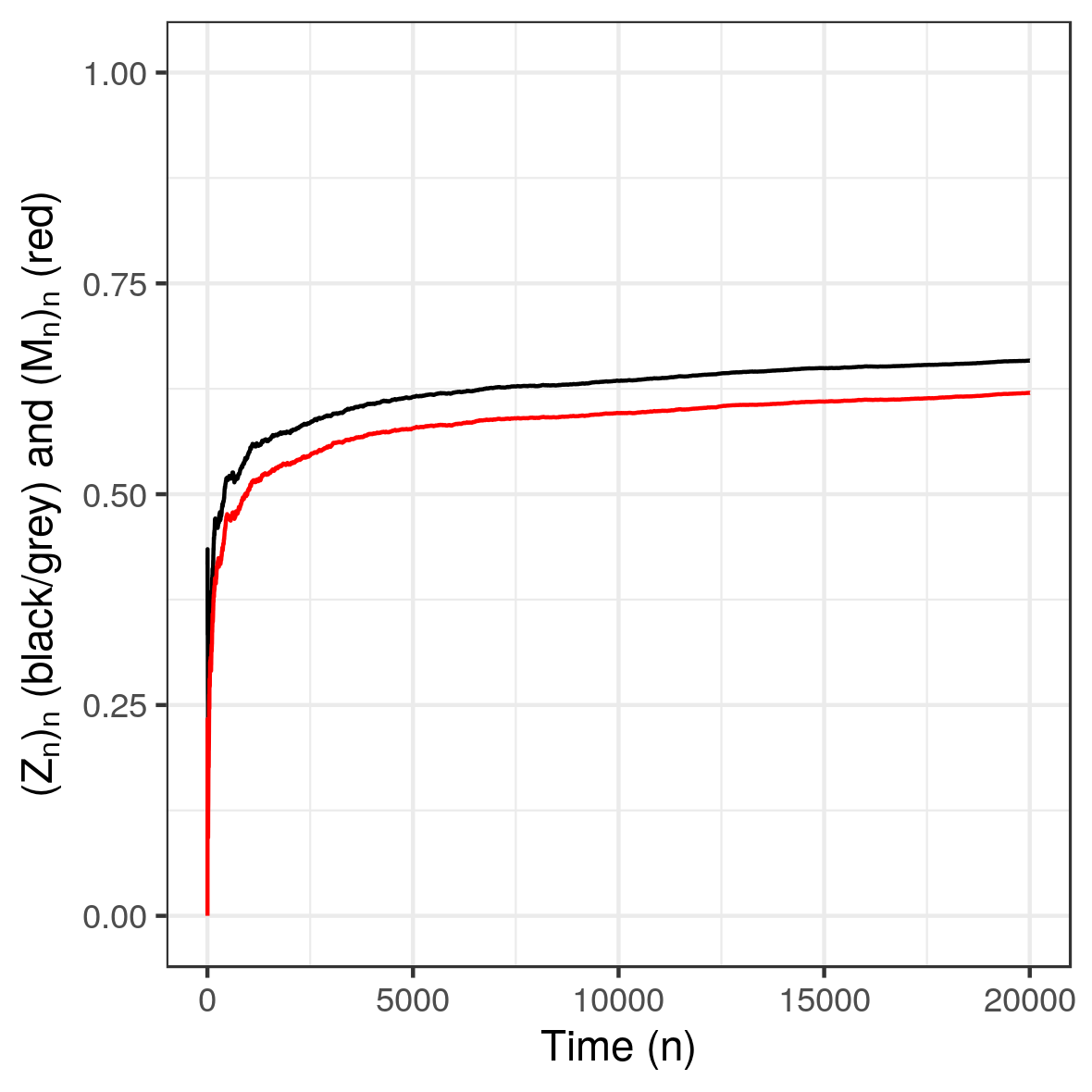}
\end{subfigure}
\caption{Case 3a. Time-horizon $5.000$ (left), $20.000$ (right).  On
  each picture, one sample plot of $(Z_n)_n$ and $(M_n)_n$ (resp. black and red).
%
  %
}
\label{fig:example4a}
\end{figure}

\noindent \textbf{Example 3a}\\
%
%
Take each $N_n$ independent of
$\mathcal{F}_{n-1}$ and uniformly distributed on
$\{1,\hdots,5\}$. Take $[A_n,B_n]$ satisfying assumption (A3) and
taking values $(1,1),\,(3,1),\,(1,3),\,(3,3)$ with respective
probabilities $\frac{3}{16},\frac14,\frac{1}{16},\frac12$.  It holds
$m_A=2.5$ and $m_B=2.125$. We set $a=b=5$.  See
Fig.~\ref{fig:example4a} for samples.  \\

\begin{figure}
\centering
\begin{subfigure}[t]{0.49\linewidth}
\centering
\includegraphics[scale=0.85]{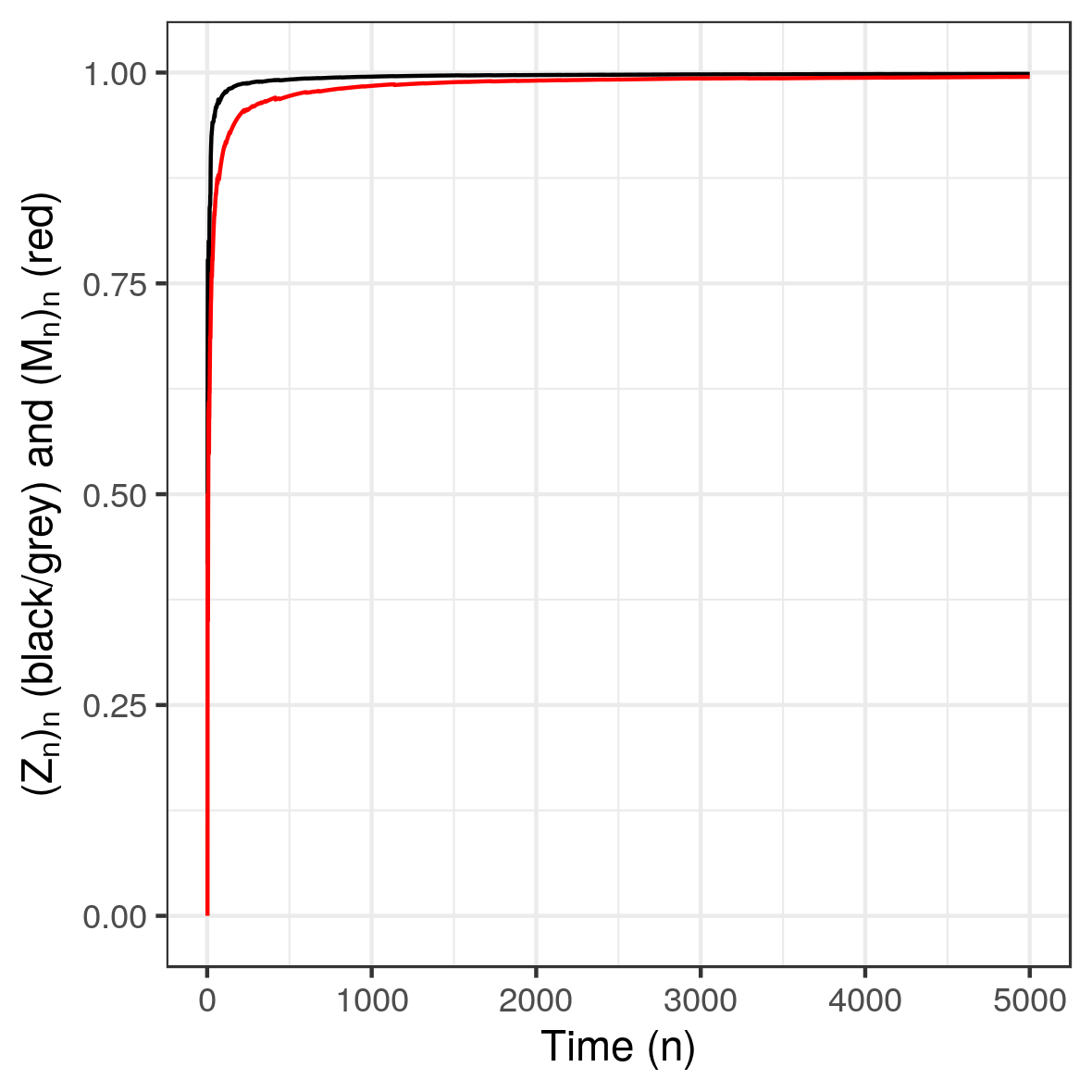}
\end{subfigure}
\hfill 
\begin{subfigure}[t]{0.49\linewidth}
\centering
\includegraphics[scale=0.85]{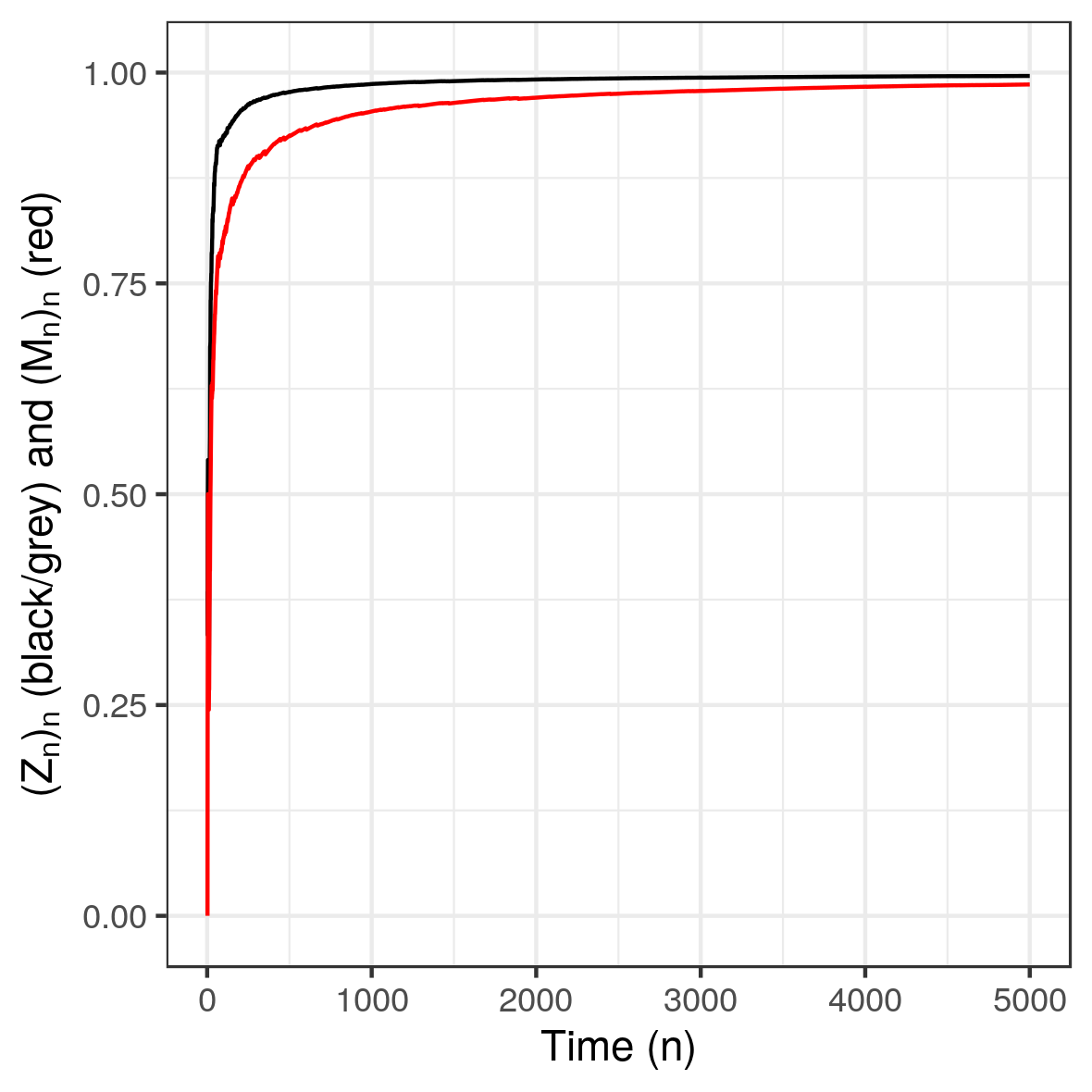}
\end{subfigure}
\caption{Case 3b. Time-horizon $5.000$.  On each picture, one sample
  plot of $(Z_n)_n$ and $(M_n)_n$ (resp. black and red).
%
}
\label{fig:example4b}
 \end{figure}

 \noindent \textbf{Example 3b} \\
 %
 %
 Take each $N_n$ independent of
 $\mathcal{F}_{n-1}$ and uniformly distributed on $\{1,\hdots,5\}$.
 Take $[A_n,B_n]$ satisfying assumption~(A3) and taking values
 $(1,1),\,(10,1),\,(1,3),\,(10,3)$ with respective probabilities
 $\frac{1}{5},\frac{2}{5},\frac{1}{5},\frac15$.  It holds $m_A=6.4$
 and $m_B=1.8$. We set $a=b=5$.  See Fig.~\ref{fig:example4b} for
 samples.  \\
 
\clearpage

\noindent{\bf Acknowledgments}\\ Irene Crimaldi and Ida Minelli are
members of the Italian Group ``Gruppo Nazionale per l'Analisi
Matematica, la Probabilit\`a e le loro Applicazioni'' of the Italian
Institute ``Istituto Nazionale di Alta Matematica''. P.-Y. Louis
acknowledges the International Associated Laboratory Ypatia Laboratory
of Mathematical Sciences (LYSM) for funding travel
expenses.
\\

\noindent{\bf Funding Sources}\\ Irene Crimaldi is partially supported
by the Italian ``Programma di Attivit\`a Integrata'' (PAI), project ``TOol
for Fighting FakEs'' (TOFFE) funded by IMT School for Advanced Studies
Lucca.
\\

\noindent{\bf Declaration}\\ All the authors equally contributed to
this work.


\appendix

\section{Technical results}

Consider the model and the assumptions described in Section~\ref{model}.

\begin{lemma}\label{lemma-ipergeom}
  Suppose $A_n\vee B_n\vee N_n\leq C$ for some (integer) constant
  $C$. Let $p_{n+1,k}=p_k(N_{n+1},S_n,H_n)$ be the values of the
  hypergeometric distribution with parameters $N_{n+1},\,S_n$ and
  $H_n$ (see \eqref{hypergeometric-distr}). Then, we have
  $$1-p_{n+1,N_{n+1}}=\frac{K_n}{S_n}\left(1+O(1)\right)=O(K_n/S_n)\,.$$
\end{lemma}

\begin{proof}
  If $N_{n+1}=1$, we simply have $1-p_{n+1,N_{n+1}}=K_n/S_n$. By Lemma
 ~\ref{colori-a-infinito}, we have $H_n\geq C$ for $n$ large enough
  (and so $H_n\geq N_{n+1}$ for $n$ large enough). Therefore, for $n$
  large enough, we have
  \begin{equation*}
    \begin{split}
      1-p_{n+1,N_{n+1}}&=1-\prod_{j=1}^{N_{n+1}}\frac{H_n-j+1}{S_n-j+1}=
      \frac{H_n+K_n}{S_n}- \prod_{j=1}^{N_{n+1}}\frac{H_n-j+1}{S_n-j+1}\\
      &=\frac{K_n}{S_n}+
      \frac{H_n}{S_n}\left[1-\prod_{j=2}^{N_{n+1}}\frac{H_n-j+1}{S_n-j+1} \right]\\
      &=\frac{K_n}{S_n}+
      \frac{H_n}{S_n}
      \frac{\prod_{j=1}^{N_{n+1}-1}(S_n-j)-\prod_{j=1}^{N_{n+1}-1}(H_n-j)}
           {\prod_{j=1}^{N_{n+1}-1}(S_n-j)}\\
&=\frac{K_n}{S_n}+
      \frac{H_n}{S_n}
      \frac{(S_n-H_n)f(S_n,H_n)}{\prod_{j=1}^{N_{n+1}-1}(S_n-j)} 
      =\frac{K_n}{S_n}
      \left(1+\frac{H_nf(S_n,H_n)}{\prod_{j=1}^{N_{n+1}-1}(S_n-j)}\right)\,, 
  \end{split}
  \end{equation*}
  where $f(x,y)=1$ when $N_{n+1}=2$ and
  $f(x,y)=\sum_{j=1}^{N_{n+1}-2}a_jx^j+b_jy^j+c$ when $N_{n+1}\geq
  3$. Therefore, since $H_n\leq S_n$ and $S_n\to +\infty$ almost
  surely (by Lemma~\ref{colori-a-infinito}), we have
  $H_nf(S_n,H_n)/\prod_{j=1}^{N_{n+1}-1}(S_n-j)=O(1)$.
\end{proof}

\begin{lemma}\label{supermart}
  Suppose to be in case 2).  For $e>1$,
  $H_{n}/K_n^e$ and $K_n/H_n^{e}$ are eventually (positive)
  supermartingales and so they converge almost surely to a finite
  random variable.
  \end{lemma}
\begin{proof}
  The proof used in order to prove that $Q_n=K_n/H_n^e$ is eventually
  a positive supermartingale in the proof of Theorem
 ~\ref{as-conv-medie-diverse} does not work now, because we have $e>1$
  and the inequality $(1-x)^e\leq 1-ex$ is not true. Therefore we need
  a different proof. We observe that
  \begin{equation*}
    \begin{split}
      &E\left[\frac{H_{n+1}}{K_{n+1}^e}-\frac{H_n}{K_n^e}\,|\,\mathcal{H}_n\right]=
      E\left[\frac{H_{n+1}}{K_n^e}-\frac{H_n}{K_n^e} +
        \frac{H_{n+1}}{K_{n+1}^e}-\frac{H_{n+1}}{K_n^e}\,|\,\mathcal{H}_n\right]=\\
      &\sum_{k\in\mathcal{X}_{n+1}}
      p_{n+1,k}\left(\frac{H_n+A_{n+1}k}{K_n^e}-\frac{H_n}{K_n^e}\right)
      +p_{n+1,k}(H_n+A_{n+1}k)
      \left(\frac{1}{(K_n+B_{n+1}(N_{n+1}-k))^e}-\frac{1}{K_n^e}\right)=\\
      &
\sum_{k\in \mathcal{X}_{n+1}} p_{n+1,k}\frac{A_{n+1}k}{K_n^e}
+
\sum_{k\in \mathcal{X}_{n+1}\setminus \{N_{n+1}\}} p_{n+1,k}(H_n+A_{n+1}k)
      \left(\frac{1}{(K_n+B_{n+1}(N_{n+1}-k))^e}-\frac{1}{K_n^e}\right)
      \,.
      \end{split}
  \end{equation*}
  Using the Taylor expansion of the function $f(x)=1/(c+x)^e$ with
  $c=K_n$ and $x=B_{n+1}(N_{n+1}-k)$, we can choose a constant
  $\theta$ such that eventually
  $$
  \left(\frac{1}{(K_n+B_{n+1}(N_{n+1}-k))}-\frac{1}{K_n^e}\right)
  \leq
  -\frac{e}{K_n^{e+1}}\left(B_{n+1}(N_{n+1}-k)-\frac{\theta}{K_n}\right)\,.
  $$ Therefore the last term of the above equalities is eventually
  smaller or equal than
  $$
  \frac{H_n}{K_n^e}\left\{
  \sum_{k\in \mathcal{X}_{n+1}}
 \left( \frac{A_{n+1}k}{H_n}- e\frac{B_{n+1}(N_{n+1}-k)}{K_n} \right)p_{n+1,k}
 + e\theta \sum_{k\in \mathcal{X}_{n+1}\setminus\{N_{n+1}\}}
 \frac{(1+A_{n+1}k/H_n)}{K_n^2}p_{n+1,k}\right\}\,.
$$
Now, we observe that
  $$
  E\left[\sum_{k\in \mathcal{X}_{n+1}}
    \left(\frac{A_{n+1}k}{H_n}-e \frac{B_{n+1}(N_{n+1}-k)}{K_n}\right)p_{n+1,k}
    \,|\,\mathcal{G}_n\right]
=m_{n+1}\frac{N_{n+1}}{S_n}(1-e)
$$
and (using $\lim_n m_n=m>0$, $N_{n+1}\geq 1$ and Lemma~\ref{lemma-ipergeom})
$$
E\left[\sum_{k\in \mathcal{X}_{n+1}\setminus\{N_{n+1}\}}
  \frac{(1+A_{n+1}k/H_n)}{K_n^2}p_{n+1,k}\,|\,\mathcal{G}_n\right]
\leq \frac{(1-p_{n+1,N_{n+1}})}{K_n^2}+\frac{m_{n+1}N_{n+1}}{S_n K_n^2}
=O(1/(S_n K_n))\,.
$$
Therefore, we have 
$$
E\left[\frac{H_{n+1}}{K_{n+1}^e}-\frac{H_n}{K_n^e}\,|\,\mathcal{G}_n\right]
\leq
m_{n+1}\frac{H_n}{K_n^e}\frac{N_{n+1}}{S_n}
\left[-(e-1)+O(1/K_n)\right]
$$ and so, since $e>1$ and $K_n\uparrow +\infty$ (by Lemma
~\ref{colori-a-infinito}), we can conclude that the above conditional
expectation is definitely negative.
%
%
\end{proof}

\begin{lemma}\label{rate-colori-2}
Under the assumptions of Theorem~\ref{th-legge-1}, we have
$1/K_n=O(1/n^\gamma)$ and $1/H_n=O(1/n^\gamma)$ for some $\gamma>0$.
\end{lemma}
\begin{proof} This proof is essentially the same as the one of Lemma A.1(iv)
 in~\cite{mayflournoy}. However, for the reader's
 convenience, we here rewrite it with all the details.  Since
 $S_n/n=(H_n+K_n)/n$ converges almost surely to $mN$, we have that
 eventually $S_n=(H_n+K_n)>n m N 3/4$ almost surely. Let $F_H=\{H_n>n
 m N/4\;\mbox{ eventually}\}$ and $F_K=\{K_n> n m N/4\;
 \mbox{eventually}\}$.  Since $(Z_n)$ converges almost surely to $Z$
 with values in $[0,1]$, then $H_n/K_n=Z_n/(1-Z_n)$ converges almost
 surely to a random variable with values in $[0,+\infty]$. It follows
 that $P(F_H\cup F_K)=1$. Indeed, on $(F_H\cup F_K)^c=F_H^c\cap
 F_K^c$, we have $\liminf H_n/n\leq mN/4$, $\liminf_n K_n/n\leq mN/4$
 and $H_n+K_n>n m N 3/4$ almost surely and so, since we can write
 $K_n/H_n=(H_n+K_n)/H_n-1$ and $H_n/K_n=(H_n+K_n)/K_n-1$, we have
 $\liminf_n H_n/K_n\leq 1/2<2\leq \limsup_n H_n/K_n$. This means that
 on $(F_H\cup F_K)^c$, $H_n/K_n$ does not converge and hence
 $P((F_H\cap F_K)^c)=0$. In order to conclude, it is enough to prove
 that on $F_H$ (resp. $F_K$), $K_n$ (resp. $H_n$) is eventually
 greater than $n^\gamma$ for $\gamma>0$ (up to a multiplicative
 constant).  \\ \indent Now, by Lemma~\ref{supermart}, $H_n/K_n^e$ is
 bounded and we know that $K_n\uparrow +\infty$ (see Lemma
~\ref{colori-a-infinito}). Therefore, for each $\epsilon>0$, we have
 $H_n/K_n^{e+\epsilon}\to 0$ almost surely and so
 $H_n/K_n^{e+\epsilon}<1$ eventually. Therefore on $F_H$, we
 eventually have $K_n^{e+\epsilon}=(H_n/K_n^{e+\epsilon})^{-1} H_n>n m
 N /4\geq n m /4$, \textit{i.e.} $K_n>n^{\gamma}$ eventually (up to a
 multiplicative constant) with $\gamma=1/(e+\epsilon)>0$. Similarly,
 on $F_K$, we have $H_n>n^{\gamma}$ eventually (up to a multiplicative
 constant) with $\gamma=1/(e+\epsilon)>0$.
\end{proof}


\section{Some auxiliary results}

For reader's convenience, we state here some general results: 

\begin{lemma}\label{lemma-app} (Lemma 2 in~\cite{BCPR-urne1})\\
Let $(Y_n)$ be a sequence of real random variables, adapted to a
filtration $\mathcal F$. If $\sum_{j\geq 1} j^{-2} E[Y_j^2]<+\infty$
and $E[Y_{j}| {\mathcal F}_{j-1}]\stackrel{a.s.}\longrightarrow Y$ for
some real random variable $Y$, then
$$
n\sum_{j\geq n}\frac{Y_j}{j^2}\stackrel{a.s.}\longrightarrow Y,
\qquad
\frac{1}{n}\sum_{j=1}^n Y_j\stackrel{a.s.}\longrightarrow Y.
$$
\end{lemma}

\begin{lemma}\label{lemma-app-C} (Th.~2 in~\cite{Bl-Du} or  
a special case of Lemma A.2 in~\cite{crimaldi-2009})\\
Let $\mathcal F$ be a filtration and set ${\mathcal
  F}_{\infty}=\bigvee_n {\mathcal F}_n$.  Then, for each sequence
$(Y_n)$ of integrable complex random variables, which is dominated in
$L^1$ and which converges almost surely to a complex random variable
$Y$, the conditional expectation $E[Y_n|{\mathcal F}_n]$ converges
almost surely to the conditional expectation $E[Y| {\mathcal
    F}_{\infty} ]$.
\end{lemma}

\section{Stable convergence and its variants}
This brief appendix contains some basic definitions and results
concerning stable convergence and its variants. For more details, we
refer the reader to \cite{crimaldi-2009, crimaldi-libro,
  cri-let-pra-2007, hall-1980} and the references therein.\\

\indent Let $(\Omega, {\mathcal A}, P)$ be a probability space,
and let $S$ be a Polish space, endowed with its Borel
$\sigma$-field. A {\em
  kernel} on $S$, or a random probability measure on $S$, is a
collection $K=\{K(\omega):\, \omega\in\Omega\}$ of probability
measures on the Borel $\sigma$-field of $S$ such that, for each
bounded Borel real function $f$ on $S$, the map
$$
\omega\mapsto K\!f(\omega)=\int f (x)\, K(\omega)(dx)
$$
is $\mathcal A$-measurable. Given a sub-$\sigma$-field $\mathcal
H$ of $\mathcal A$, a kernel $K$ is said $\mathcal H$-measurable
if all the
above random variables $K\!f$ are $\mathcal H$-measurable.\\

\indent On $(\Omega, {\mathcal A},P)$, let $(Y_n)_n$ be a sequence
of $S$-valued random variables, let $\mathcal H$ be a
sub-$\sigma$-field of $\mathcal A$, and let $K$ be a $\mathcal
H$-measurable kernel on $S$. Then we say that $Y_n$ converges {\em
$\mathcal H$-stably} to $K$, and we write $Y_n\longrightarrow K$
${\mathcal H}$-stably, if
$$
P(Y_n \in \cdot \,|\, H)\stackrel{weakly}\longrightarrow
E\left[K(\cdot)\,|\, H \right] \qquad\hbox{for all } H\in{\mathcal
H}\; \hbox{with } P(H) > 0,
$$where $K(\cdot)$ denotes the random variable
  defined, for each Borel set $B$ of $S$, as $\omega\mapsto
  K\!I_B(\omega)=K(\omega)(B)$.  In the case when ${\mathcal
  H}={\mathcal A}$, we simply say that $Y_n$ converges {\em stably} to
$K$ and we write $Y_n\longrightarrow K$ stably. Clearly, if
$Y_n\longrightarrow K$ ${\mathcal H}$-stably, then $Y_n$ converges
in distribution to the probability distribution $E[K(\cdot)]$.
Moreover, the $\mathcal H$-stable convergence of $Y_n$ to $K$ can
be stated in terms of the following convergence of conditional
expectations:
\begin{equation}\label{def-stable}
E[f(Y_n)\,|\, {\mathcal H}]\stackrel{\sigma(L^1,\,
L^{\infty})}\longrightarrow K\!f
\end{equation}
for each bounded continuous real function $f$ on $S$. \\

\indent in~\cite{cri-let-pra-2007} the notion of $\mathcal
H$-stable convergence is firstly generalized in a natural way
replacing in (~\ref{def-stable}) the single sub-$\sigma$-field
$\mathcal H$ by a collection ${\mathcal G}=({\mathcal G}_n)_n$
(called conditioning system) of sub-$\sigma$-fields of $\mathcal
A$ and then it is strengthened by substituting the convergence in
$\sigma(L^1,L^{\infty})$ by the one in probability (\textit{i.e.} in $L^1$,
since $f$ is bounded). Hence, according to
\cite{cri-let-pra-2007}, we say that $Y_n$ converges to $K$ {\em
stably in the strong sense}, with respect to ${\mathcal
G}=({\mathcal G}_n)_n$, if
\begin{equation}\label{def-stable-strong}
E\left[f(Y_n)\,|\,{\mathcal G}_n\right]\stackrel{P}\longrightarrow
K\!f
\end{equation}
for each bounded continuous real function $f$ on $S$.\\

\indent Finally, a strengthening of the stable convergence in the
strong sense can be naturally obtained if in (\ref{def-stable-strong})
we replace the convergence in probability by the almost sure
convergence (see \cite{crimaldi-2009}): given a conditioning system
${\mathcal G}=({\mathcal G}_n)_n$, we say that $Y_n$ converges to $K$
in the sense of the {\em almost sure conditional convergence}, with
respect to ${\mathcal G}$, if
\begin{equation*}
E\left[f(Y_n)\,|\,{\mathcal
G}_n\right]\stackrel{a.s.}\longrightarrow K\!f
\end{equation*}
for each bounded continuous real function $f$ on
  $S$. 
\\

\indent We conclude recalling two results. In particular, for the
second one, we denote by $\mathcal{N}(\mu,\sigma^2)$ the Gaussian
probability distribution with mean $\mu$ and variance $\sigma^2\geq 0$
(where $\mathcal{N}(\mu, 0)$ means the Dirac distribution $\delta_\mu$
concentrated in $\mu$). Therefore, when $U$ is a positive random
variable, the symbol $\mathcal{N}(0,U)$ denotes the Gaussian kernel
$\{\mathcal{N}(0,U(\omega)):\, \omega\in\Omega\}\}$.

\begin{theorem}\label{lemma-conv-stab}
  (Lemma 1 in~\cite{BCPR-urne1})\\ Suppose that $C_n$ and $D_n$ are
  $S$-valued random variables, that $M$ and $N$ are kernels on $S$,
  and that $\mathcal{G}=(\mathcal{G}_n)_n$ is a filtration satisfying
  $\sigma(C_n)\underline\subset\mathcal{G}_n$ and
  $\sigma(D_n)\underline\subset\sigma( \cup_n\mathcal{G}_n)$ for all
  $n$. If $C_n$ stably converges to $M$ and $D_n$ converges to $N$
  stably in the strong sense, with respect to $\mathcal{G}$, then
  $[C_n,D_n ]\longrightarrow M \otimes N$ stably. (Here, $M
  \otimes N$ is the kernel on $S \times S$ such that $(M \otimes N
  )(\omega) = M(\omega) \otimes N (\omega)$ for all $\omega$.)
\end{theorem}
This last result contains as a special case the fact that stable
convergence and convergence in probability combine well: that is, if
$C_n$ stably converges to $M$ and $D_n$ converges in probability to a
random variable $D$, then $(C_n, D_n)$ stably converges to $M \otimes
\delta_D$, where $\delta_D$ denotes the Dirac kernel concentrated in
$D$. In particular, if $M$ is the Gaussian kernel $\mathcal{N}(0,D)$,
we have $C_n/\sqrt{D_n}\longrightarrow \mathcal{N}(0,1)$ stably.

\begin{theorem}\label{th-main-app} 
  (See Th.~1 together with Prop.~1 in~\cite{BCPR-urne1} and Th.~10 in
  \cite{BCPR-indian})\\ Let $(Y_n)$ be a bounded sequence of real
  random variables, adapted to a filtration ${\mathcal
    G}=(\mathcal{G}_n)$. Set
\begin{equation*}
Z_n=E[Y_{n+1}|\mathcal{G}_n]\quad\text{and}\quad 
M_n=\frac{1}{n}\sum_{j=1}^n Y_j.
\end{equation*}
Suppose that $n^3 E\left[\,(E[Z_{n+1}|\mathcal{G}_n]-Z_n)^2\,\right]\to 0$.
\\
\indent Then, $Z_n\overset{a.s.}\longrightarrow Z$ and
$M_n\overset{a.s.}\longrightarrow Z$ for some real random variable
$Z$. Moreover, $\sqrt{n}(Z_n-Z)$ converges in the sense of the almost
sure conditional convergence with respect to $\mathcal G$ toward the
Gaussian kernel~${\mathcal N}(0,V)$ for some real random variable $V$,
provided
\begin{itemize}
\item[c1)] $E\left[\sup_{j\geq
    1}\sqrt{j}\,|Z_{j-1}-Z_j|\,\right]<+\infty$, 
\item[c2)] $n\sum_{j\geq n}(Z_{j-1}-Z_j)^2\stackrel{a.s.}\longrightarrow
V$.
\end{itemize}
If condition
\begin{itemize}
\item[c3)] $n^{-1}\sum_{j=1}^n\bigl[Y_j-Z_{j-1}+j(Z_{j-1}-Z_j)\bigr]^2
  \stackrel{P}\longrightarrow U$
\end{itemize}
is also satisfied for some real random variable $U$, then
\begin{equation*}
\left[\sqrt{n}\,\bigl(M_n-Z_n\bigr),\sqrt{n}(Z_n-Z)\right]
\stackrel{stably}\longrightarrow
\mathcal{N}\bigl(0,\,U\bigr)\otimes{\mathcal N}(0,V).
\end{equation*}
\end{theorem}
In particular, we have $\sqrt{n}\,\bigl(M_n-Z_n\bigr)\longrightarrow
\mathcal{N}(0,U)$ stably and
$\sqrt{n}\,\bigl(M_n-Z\bigr)\longrightarrow \mathcal{N}(0,U+V)$
stably.

%
%


\begin{thebibliography}{10}

\bibitem{AguechLasmarSelmi2019}
R.~Aguech, N.~Lasmar, and O.~Selmi.
\newblock A generalized urn with multiple drawing and random addition.
\newblock {\em Annals of the Institute of Statistical Mathematics},
  71(2):389--408, Apr. 2019.

\bibitem{AguechSelmi-unbalanced}
R.~Aguech and O.~Selmi.
\newblock Unbalanced multi-drawing urn with random addition matrix.
\newblock {\em Arab Journal of Mathematical Sciences}, 2019.

\bibitem{perron}
D.~A. Aoudia and F.~Perron.
\newblock A new randomized {P}\'olya urn model.
\newblock {\em Appl. Math.}, 3:2118--2122, 2012.

\bibitem{Berti2010}
P.~Berti, I.~Crimaldi, L.~Pratelli, and P.~Rigo.
\newblock {Central limit theorems for multicolor urns with dominated colors}.
\newblock {\em Stoch. Process. their Appl.}, 120(8):1473--1491, 2010.

\bibitem{BCPR-urne1}
P.~Berti, I.~Crimaldi, L.~Pratelli, and P.~Rigo.
\newblock A central limit theorem and its applications to multicolor randomly
  reinforced urns.
\newblock {\em Journal of Applied Probability}, 48(2):527--546, 2011.

\bibitem{BCPR-indian}
P.~Berti, I.~Crimaldi, L.~Pratelli, and P.~Rigo.
\newblock Central limit theorems for an {I}ndian buffet model with random
  weights.
\newblock {\em The Annals of Applied Probability}, 25(2):523--547, 2015.

\bibitem{Bl-Du}
D.~Blackwell and L.~Dubins.
\newblock Merging of opinions with increasing information.
\newblock {\em The Annals of Mathematical Statistics}, 33(3):882--886, 1962.

\bibitem{Chen2020}
M.-R. Chen.
\newblock A time dependent {P}{\'{o}}lya urn with multiple drawings.
\newblock {\em Probab. Eng. Informational Sci.}, 34(4):469--483, 2020.

\bibitem{chen-kuba}
M.-R. Chen and M.~Kuba.
\newblock On generalized {P}\'olya urn models.
\newblock {\em J. Appl. Prob.}, 50:1169--1186, 2013.

\bibitem{chen-wei-2005}
M.-R. Chen and C.-Z. Wei.
\newblock A new urn model.
\newblock {\em Journal of Applied Probability}, 42(4):964--976, Dec. 2005.

\bibitem{crimaldi-2009}
I.~Crimaldi.
\newblock An almost sure conditional convergence result and an application to a
  generalized {P}{\'o}lya urn.
\newblock {\em Int. Math. Forum}, 4(21-24):1139--1156, 2009.

\bibitem{crimaldi2016}
I.~Crimaldi.
\newblock Central limit theorems for a hypergeometric randomly reinforced urn.
\newblock {\em Journal of Applied Probability}, 53(3):899--913, 2016.

\bibitem{crimaldi-libro}
I.~Crimaldi.
\newblock {\em {I}ntroduzione alla nozione di convergenza stabile e sue
  varianti ({I}ntroduction to the notion of stable convergence and its
  variants)}, volume~57.
\newblock {U}nione {M}atematica {I}taliana, Monograf s.r.l., Bologna, Italy.,
  2016.
\newblock Book written in Italian.

\bibitem{cri-dai-min}
I.~Crimaldi, P.~Dai~Pra, and I.~G. Minelli.
\newblock Fluctuation theorems for synchronization of interacting {P}\'olya's
  urns.
\newblock {\em Stochastic Process. Appl.}, 126(3):930--947, 2016.

\bibitem{Crimaldi2008}
I.~Crimaldi and F.~Leisen.
\newblock {Asymptotic Results for a Generalized P{\'{o}}lya Urn with
  "Multi-Updating" and Applications to Clinical Trials}.
\newblock {\em Commun. Stat. - Theory Methods}, 37(17):2777--2794, July 2008.

\bibitem{cri-let-pra-2007}
I.~Crimaldi, G.~Letta, and L.~Pratelli.
\newblock A strong form of stable convergence.
\newblock In {\em S\'eminaire de {P}robabilit\'es {XL}}, volume 1899 of {\em
  Lecture Notes in Math.}, pages 203--225. Springer, Berlin, 2007.

\bibitem{egg-pol}
F.~Eggenberger and G.~P\'{o}lya.
\newblock {\"{U}}ber die {S}tatistik verketteter {V}org\"{a}nge.
\newblock {\em Z. Angewandte Math. Mech.}, 3:279--289, 1923.

\bibitem{hall-1980}
P.~Hall and C.~C. Heyde.
\newblock {\em Martingale limit theory and its application}.
\newblock Academic Press, Inc. [Harcourt Brace Jovanovich, Publishers], New
  York-London, 1980.
\newblock Probability and Mathematical Statistics.

\bibitem{Higueras2006}
I.~Higueras, J.~Moler, F.~Plo, and M.~San~Miguel.
\newblock Central limit theorems for generalized {P}{\'o}lya urn models.
\newblock {\em Journal of Applied Probability}, 43(4):938--951, 2006.

\bibitem{Idriss2018}
S.~Idriss and N.~Lasmar.
\newblock {Limit Theorems for Stochastic Approximations Algorithms With
  Application to General Urn Models}.
\newblock Hal-01726014, 2018.

\bibitem{Johnson2004}
N.~Johnson, S.~Kotz, and H.~Mahmoud.
\newblock {P{\'o}lya-Type Urn Models with Multiple Drawings}.
\newblock {\em J. Iran. Stat. Soc.}, 3(2):165--173, 2004.

\bibitem{Kotz1997}
S.~Kotz and N.~Balakrishnan.
\newblock {\em {Advances in Urn Models during the Past Two Decades}},
  chapter~14, pages 203--257.
\newblock Statistics for Industry and Technology. Birkh{\"{a}}user Boston,
  1997.

\bibitem{kuba2016classification}
M.~Kuba.
\newblock Classification of urn models with multiple drawings.
\newblock Preprint Arxiv 1612.04354, 2016.

\bibitem{KubaMahmoudPanholzer}
M.~Kuba, H.~Mahmoud, and A.~Panholzer.
\newblock Analysis of a generalized {Friedman}'s urn with multiple drawings.
\newblock {\em Discrete Appl. Math.}, 161(18):2968--2984, Dec. 2013.

\bibitem{KubaMahmoud-balanced-affine-2017}
M.~Kuba and H.~M. Mahmoud.
\newblock Two-color balanced affine urn models with multiple drawings.
\newblock {\em Adv. in Appl. Math.}, 90:1--26, Sept. 2017.

\bibitem{Kuba-Sulzbach}
M.~Kuba and H.~Sulzbach.
\newblock On martingale tail sums in affine two-color urn models with multiple
  drawings.
\newblock {\em J. Appl. Probab.}, 54(1):96--117, 2017.

\bibitem{laslier2017}
B.~Laslier and J.-F. Laslier.
\newblock Reinforcement learning from comparisons: {T}hree alternatives are
  enough, two are not.
\newblock {\em Ann. Appl. Probab.}, 27(5):2907--2925, Oct. 2017.

\bibitem{mailler}
N.~Lasmar, C.~Mailler, and O.~Selmi.
\newblock {Multiple drawing multi-colour urns by stochastic approximation}.
\newblock {\em J. Appl. Probab.}, 55(1):254--281, 2018.

\bibitem{launay2012urns}
M.~Launay.
\newblock Urns with simultaneous drawing.
\newblock Preprint Arxiv 1201.3495, 2012.

\bibitem{mah}
H.~M. Mahmoud.
\newblock {\em {P}\'olya urn models}.
\newblock Texts in Statistical Science Series. CRC Press, Boca Raton, FL, 2009.

\bibitem{mahmoud_2013_multisets}
H.~M. Mahmoud.
\newblock Drawing multisets of balls from tenable balanced linear urns.
\newblock {\em Probability in the Engineering and Informational Sciences},
  27(2):147--162, 2013.

\bibitem{mayflournoy}
C.~May and N.~Flournoy.
\newblock Asymptotics in response-adaptive designs generated by a two-color,
  randomly reinforced urn.
\newblock {\em Ann. Statist.}, 37(2):1058--1078, 04 2009.

\bibitem{EpidemicNetworks}
R.~Pastor-Satorras, C.~Castellano, P.~Van~Mieghem, and A.~Vespignani.
\newblock Epidemic processes in complex networks.
\newblock {\em Rev. Modern Phys.}, 87(3):925--979, 2015.

\bibitem{pem}
R.~Pemantle.
\newblock A survey of random processes with reinforcement.
\newblock {\em Probab. Surv.}, 4(1-79):1--79, 2007.

\bibitem{pemantle-volkov-1999}
R.~Pemantle and S.~Volkov.
\newblock Vertex-reinforced random walk on $\mathbf{Z}$ has finite range.
\newblock {\em Ann. Probab.}, 27(3):1368--1388, July 1999.

\end{thebibliography}

\end{document}